\documentclass[10pt]{article}
\usepackage{amsfonts,amsmath,amsthm,amssymb,tikz}
\usepackage{pgfplots}
\usepackage{graphics, epsfig}
\usepackage{color}
\usepackage{appendix}
\usepackage{verbatim}
\usepackage{ulem}
\usepackage[makeroom]{cancel}
\usepackage[margin=1.5in]{geometry}
\usepackage[colorlinks=true, urlcolor=blue, citecolor=red, linkcolor=blue, pdffitwindow=true, linktocpage, pdfpagelabels, bookmarksnumbered, bookmarksopen]{hyperref}

 \usepackage[usenames,dvipsnames]{pstricks}
 \usepackage{pst-grad} % For gradients
 \usepackage{pst-plot} % For axes
\allowdisplaybreaks

\newtheorem{remark}{Remark}[section]

\def\Xint#1{\mathchoice
{\XXint\displaystyle\textstyle{#1}}%
{\XXint\textstyle\scriptstyle{#1}}%
{\XXint\scriptstyle\scriptscriptstyle{#1}}%
{\XXint\scriptscriptstyle\scriptscriptstyle{#1}}%
\!\int}
\def\XXint#1#2#3{{\setbox0=\hbox{$#1{#2#3}{\int}$ }
\vcenter{\hbox{$#2#3$ }}\kern-.6\wd0}}

\def\dashint{\Xint-}

\newtheorem{theorem}{Theorem}[section]
\newtheorem{lemma}{Lemma}[section]

\numberwithin{equation}{section}
\numberwithin{theorem}{section}
\numberwithin{proposition}{section}
\numberwithin{lemma}{section}
\numberwithin{remark}{section}
\newcommand{\osc}{\operatornamewithlimits{osc}}

\def\R{\mathbb{R}}
\def\Q{\mathcal{Q}}
\def\K{\mathcal{K}}
\def\N{\mathbb{N}}

\def\d{\mathrm{d}}

\def\A{\mathcal{A}}
\def\S{\mathcal{S}}

\def\e{\epsilon}
\begin{document}
\title{On the Local Behavior of Local Weak Solutions to some Singular Anisotropic Elliptic Equations}
%%%%%%%%%%%%%%%%%%%%%%%%%%%%%%%%%
%%%%%%%%%%%%%%%%%%%%%%%%%%%%%%%%%%%%%%%%%%%%%%%%%%%
\author{
%Kateryna O. Buryachenko,\\
%Vasyl Stus Donetsk National University,\\ 600-richa Str., 21, Vinnytsia, 21021, Ukraine\\  $katarzyna \textunderscore @ukr.net$
\\
\\
\it{Simone Ciani \& Igor I. Skrypnik \& Vincenzo} Vespri \\ \\ Universit\`a degli Studi di Firenze,\\ Dipartimento di Matematica e Informatica "Ulisse Dini"\\   $simone.ciani@unifi.it$ \& $vincenzo.vespri@unifi.it$\\ \\
Institute of Applied Mathematics and Mechanics,\\
National Academy of Sciences of Ukraine,\\
Gen. Batiouk Str. 19, 84116 Sloviansk, Ukraine\\
$iskrypnik@iamm.donbass.com$
}
%%%%%%%%%%%%%%%%%%%%%%%%%%%%%%%%%%%%%%%%%%%%%%%%%%%
%%%%%%%%%%%%%%%%%%%%%%%%%%%%%%%%%%%%
\date{}
\maketitle
\vskip.4truecm
%%%%%%%%%%%%%%%%%%%%%%%%%%%%%%%%%%%%%%%%%%%%%%%%%%%
%%%%%%%%%%%%%%%%%%%%%%%%%%%%%%%%%%%%%%%%%%%%
\begin{abstract} \noindent
We study the local behavior of bounded local weak solutions to a class of anisotropic singular equations of the kind
\begin{equation*} \label{prototype}
\sum_{i=1}^s \partial_{ii} u+ \sum_{i=s+1}^N \partial_i \bigg(A_i(x,u,\nabla u) \bigg)  =0,\quad x \in \Omega \subset \subset \R^N \quad \text{ for } \quad 1\leq s \leq (N-1),
\end{equation*}
where each operator $A_i$ behaves directionally as the singular $p$-Laplacian, $1< p < 2$. Throughout a parabolic approach to expansion of positivity we obtain the interior H\"older continuity, and some integral and pointwise Harnack inequalities. 
%%%%%%%%%%%%%%%%%%%%%%%%%%%%%%%%%%%%%%%%%%%%%%%%%%%
\vskip.2truecm
%%%%%%%%%%%%%%%%%%%%%%%%%%%%%%%%%%%%%%%%%%%%%%%%%%%
\noindent {{\bf MSC 2020:}} 
35J75, 35K92, 35B65.
\vskip.2truecm
%%%%%%%%%%%%%%%%%%%%%%%%%%%%%%%%%%%%%%%%%%%%%%%%%%%
\noindent {{\bf Key Words}}: 
Anisotropic $p$-Laplacian, Singular Parabolic Equations, H\"older Continuity, Intrinsic Scaling, Expansion of Positivity, Intrinsic Harnack Inequality.
%%%%%%%%%%%%%%%%%%%%%%%%%%%%%%%%%%%%%%%%%%%%%%%%%%%
\\
% \begin{flushright}
% \it{To celebrate Umberto Mosco's 80th genethliac}
% \end{flushright}
\end{abstract}

\newpage
%%%%%%%%%%%%%%%%%%%%%%%%%%%%%%%%%%%%%%%%%%%%%%%%%%%%%%%%%%%%%%%%%%%%%%%
 \begin{center} \section{Introduction}\label{S:intro}  \end{center}
In this note we study local regularity properties for bounded weak local solutions to operators whose prototype is 
\begin{equation}\label{prototipo}
    \sum_{i=1}^s \partial_{ii} u + \sum_{i=s+1}^N
 \partial_i \bigg( |\partial_i u|^{p-2}\partial_i u\bigg) =0, \quad \text{weakly in}\quad \Omega\subset \R^N, \quad 1<p<2, \end{equation} \noindent having a non-degenerate behavior along the first $s$- variables, and a singular behavior on the last ones. This kind of operators are useful to describe the steady states of non-Newtonian fluids that have different directional diffusions (see for instance \cite{AS}), besides their pure mathematical interest, which still is a challenge after more than fifty years. Precise hypothesis will be given later (in Section \ref{DefMain}), leaving here the space to describe what are the novelties and significance of the present work in the context of this kind of operators.\vskip0.2cm \noindent 
 Until this moment it is not known whether solutions to equations as \eqref{prototipo} enjoy the usual local properties as  $p$-Laplacean ones. This is because equation \eqref{prototipo} is part of a more general group of operators, whose regularity theory is still fragmented and largely incomplete. It is clear that new techniques are needed for a correct interpretation of the problem and its resolution. The present work is conceived to introduce a new method, adapted from the theory of singular parabolic equations. In next Section we explain this simple but effective idea, that we will apply to a class of equations as \eqref{prototipo}, that have no homogeneity on the differential operator (hence the epithet {\it anisotropic}) because they combine both non-degenerate and singular properties.\vskip0.2cm \noindent

 \subsection{The parabolic approach}
 \noindent To introduce our approach, we present an alternative proof of the Mean Value Theorem for solutions to Laplace equation. This brief and modest scheme will highlight the essence of our method, that is conceived to obtain classical properties of some elliptic equations through a parabolic approach. Let us consider the Laplace equation in an open bounded set $\Omega \subset \subset \R^N$,
 \[
 \sum_{i=1}^N \int_{\Omega} (\partial_i u)( \partial_i \phi)\, dx =0, \quad \quad \phi \in C_o^{\infty}(\Omega).
 \] Let $x_o \in \Omega$ be a Lebesgue point for $u$, and let $B_{2r}(x_o)$ be the ball of radius $2r$ and center $x_o$. Now for $0<t<r$ such that $B_{2r}(x_o) \subset \Omega$, consider the test function $\phi(t,x)=(t^2-|x-x_o|^2)_+$, to obtain the integral equality
 \[
 \sum_{i=1}^N \int_{B_t(x_o)} (\partial_i u(x)) (x_i-x_{o,i})\, dx=0.
 \] By Green's formula, this is equivalent to
%  \[
% 0= \sum_{i=1}^N \int_{B_t(x_o)} \partial_i \bigg( u\cdot  (x_i-x_o^i)\bigg) \, dx-N\int_{B_t(x_o)} u\, dx = \]
% \[
% \int_{\partial B_{t}(x_o)}  u \cdot \sum_{i=1}^N (x_i-x_o^i) \frac{(x_i-x_o^i)}{|x-x_o|} \, d\mathcal{H}^{N-1}-N\int_{B_t(x_o)} u\, dx=
% \] 
\[
\sum_{i=1}^N \int_{B_t(x_o)} \partial_i \bigg( u(x)\, (x_i-x_{o,i})\bigg) \, dx-N\int_{B_t(x_o)} u\, dx=t \int_{\partial B_t(x_o)} u\, d\mathcal{H}^{N-1}-N\int_{B_t(x_o)} u\, dx=0,\] having used that ${\bf{n}}= (x-x_o)/|x-x_o|$ is the normal unit vector to $\partial B_t(x_o)$ and being $d\mathcal{H}^{N-1}$ the Hausdorff $(N-1)$-dimensional measure. Now, last display can be rewritten as
% \[
% 0=t \int_{\partial B_t(x_o)} u\, d\mathcal{H}^{N-1}-N\int_{B_t(x_o)} u\, dx
% \]
\[
% =t^{N+1} \frac{d}{dt} \bigg(t^{-N} \int_{B_t(x_o)} u \, dx  \bigg) \quad \Rightarrow \quad 
t^{N+1}\frac{d}{dt}\bigg( t^{-N}  \int_{B_t(x_o)} u \, dx\bigg)=0.
\] Finally we integrate along $t \in (0,r)$ and we use Lebesgue's Theorem to get
\[
t^{-N}\int_{B_t(x_o)} u \, dx \bigg|_{0}^r= 
% r^{-N} \int_{B_r(x_o)} u\, dx- \lim_{t \downarrow 0} \omega_N \dashint_{B_t(x_o)}u\, dx=
r^{-N} \int_{B_r(x_o)} u\, dx-\omega_N u(x_o)=0,
\] with $\omega_N= |B_1|$ and $|B_r|= \omega_N r^{-N}$. This implies the mean value property 
\[
u(x_o) =\dashint_{B_r(x_o)} u\, dx.
\]
This point-wise control given in integral average can be used in turn to derive very strong regularity properties of the solutions. We will undergo a similar strategy for solutions to \eqref{prototipo}, by taking into account  the degeneracies and singularities that are typical of anisotropic equations.

 \subsection{Definitions and Main Results} \label{DefMain}
 \noindent 
Let $\Omega \subset \R^N$ be an open bounded set with $N \ge 2$, and let us denote with $\partial_i$ the $i$-th partial weak derivative. For $1<p<2$ and $1\leq s\leq N-1$ we consider the elliptic partial differential equation
\begin{equation}\label{E1}
    \sum_{i=1}^s \partial_{ii} u + \sum_{i=s+1}^N
 \partial_i A_i(x,u, \nabla u)=0, \quad \text{weakly in}\quad \Omega, \end{equation} where the Caratheodory\footnote{Measurable in $(u, \xi)$ for all $x \in \Omega$ and continuous in $x$ for a.e. $(u,\xi)\in \R\times \R^N$.} functions $A_i(x,u,\xi):\Omega\times  \R\times \R^N \rightarrow \R$ are subject to the following structure conditions for almost every $x \in \Omega$,
 \begin{equation}\label{E1-structure}
     \begin{cases}
      \sum_{i=s+1}^N A_i(x,u,\xi) \cdot \xi_i \ge C_1 \sum_{i=s+1}^N |\xi_i|^p -C, \quad \text{for} \quad \xi \in \R^N, \\  \\
      |A_i(x,u,\xi)| \leq C_2 |\xi_i|^{p-1}+C, \quad \text{for} \quad i\in \{ s+1,..,N\},
     \end{cases}
 \end{equation} \noindent where $C_1,C_2 >0$, $C \ge 0$ are given constants that we will always refer to as {\it the data}. A function $u \in L^{\infty}_{loc}(\Omega) \cap W^{1,[2,p]}_{loc}(\Omega)$, where
\begin{equation*}\label{anisotropic-Sobolev-spaces}
\begin{aligned}
  &W_{loc}^{1,[2,p]}(\Omega):= \bigg{\{}u\in L^1_{loc}(\Omega)\, | \, \partial_i u \in L^2_{loc}(\Omega) \,\, \forall i=1,..,s, \quad \partial_i u \in L^p_{loc}(\Omega) \,\, \forall i=s+1,..,N   \bigg{\}}, \\
  &W^{1,[2,p]}_o(\Omega):= W^{1,1}_o(\Omega) \cap   W_{loc}^{1,[2,p]}(\Omega), 
  \end{aligned}
\end{equation*}\noindent
is called a local weak solution to \eqref{E1}-\eqref{E1-structure} if for each compact set $K \subset \subset \Omega$ it satisfies
\begin{equation}\label{def-solution}
    \int \int_K \sum_{i=1}^s \partial_i u \, \partial_i \varphi \, dx + \int \int_K \sum_{s+1}^N A_i(x,u,\nabla u)\, \partial_i \varphi \, dx=0,\quad \forall \varphi \in W^{1,[2,p]}_o (K).
\end{equation}\noindent
All along the present work we will suppose that truncations $\pm (u-k)_{\pm}$ of local weak solutions to \eqref{E1}-\eqref{E1-structure} preserve the property of being sub-solutions: for any $k \in \R$, every compact subset $K \subset \Omega$, and $\psi \in W_o^{1,[2,p]}(K)$ we have
\begin{equation}\label{mon-energy}
\int \int_{K} \bigg{\{} \sum_{i=1}^s  \partial_i (u-k)_{\pm} \partial_i\psi \, + \, \sum_{i=s+1}^N A_i(x, (u-k)_{\pm}, \partial_i (u-k)_{\pm}) \partial_i \psi \, \bigg{\}} dx  \leq 0.
\end{equation} \noindent 
\begin{remark}
Previous assumption \eqref{mon-energy} is very natural. In case of homogeneous coercivity, that is, if in the first formula of \eqref{E1-structure} we have just
\[
\sum_{i=s+1}^N A_i(x,u,\xi) \cdot \xi_i \ge C_1 \sum_{i=s+1}^N |\xi_i|^p, \quad \text{for} \quad \xi \in \R^N,
\] then by a simple limit argument it can be shown that \eqref{mon-energy} is always in force (see for example \cite{DB}, Lemma 1.1 Chap. II).
\end{remark} 
\noindent Properties of anisotropic Sobolev spaces have first been investigated in \cite{KK}, \cite{Troisi},\cite{Schmeiser}, and boundedness of local weak solutions has been first considered in \cite{Kolodii} and refined in \cite{DBGV-aniso}. Limit growth conditions have been investigated in \cite{FS} and then refined in \cite{CU1}, \cite{CU2} in great generality. Henceforth it is a well-known fact in literature that local weak solutions to our equation \eqref{E1} are bounded provided $p_{max} \leq  N\bar{p}/(N-\bar{p})$, being $\bar{p}= N(\sum_{i=1}^N (p_i)^{-1})^{-1}$ the harmonic mean. We consider the prototype equation to \eqref{E1} as a special case of the full anisotropic analogue
\begin{equation}\label{full-EQ}
    -\sum_{i=1}^N \partial_i \bigg( |\partial_i u|^{p_i-2} \partial_i u \bigg)=0,
\end{equation} \noindent with $p_i=2$ for $i=1,..,s$ and $1<p_i=p<2$ on the remaining components. This last equation suffers heavily from the combined effect of singular and degenerate behavior, even when for instance all $p_i$s are greater than two. This is because the natural intrinsically scaled geometry of the equation that maintains invariant the volume $|\K|= \rho^N$ can be shaped on anisotropic cubes as
\[
\K= \prod_{i=1}^N \bigg{\{}|x_i| < M^{\frac{p_i-\bar{p}}{p_i}} \rho^{\frac{\bar{p}}{p_i}}  \bigg{\}}, 
\] where $M$ is a number depending on the solution $u$ itself (indeed the epithet {\it intrinsic}) that vanishes as soon as $u$ vanishes. Therefore when $M$ approaches zero, for those directions whose index satisfies $p_i>\bar{p}$ the anisotropic cube $\K$ shrinks to a vanishing measure, while for the remaining ones it stretches to infinity. For a detailed description of this geometry and its derivation through self-similarity we refer to \cite{Ciani}, where the evolutionary, fully anisotropic prototype equation is considered.
% We believe that a uniform theory of singular and degenerate equations where both phenomena can be interpreted in a unique general idea would solve this problem and be of great interest for a theory of anisotropic Sobolev spaces. 

\vskip0.2cm \noindent 
We state our two main results hereafter. The first one is a result of local H\"older continuity.
\begin{theorem}\label{holderTHM}
 Let $u$ be a bounded local weak solution to \eqref{E1},\eqref{E1-structure},\eqref{mon-energy}. Then there exists $\alpha\in (0,1)$ depending only on the data such that $u \in C^{0,\alpha}_{loc}(\Omega)$.
\end{theorem} \noindent 
\vskip0.2cm \noindent 
 Next we fix some geometrical notations and conventions. For a point $x_o\in \Omega$, let us denote it by  $x_o=(x_o',x_o'')$ where $x_o'\in \R^s$ and $x_o''\in \R^{N-s}$. Let $\theta,\rho>0$ be two parameters, and define the polydisc 
 \begin{equation}\label{cylinders}
     Q_{\theta,\rho}(x_o):= B_{\theta}(x_o') \times B_{\rho}(x_o'').
 \end{equation} \noindent We will say $Q_{\theta,\rho}$ is an {\it intrinsic} polydisc when $\theta$ depends on the solution $u$ itself. We will call first $s$ variables the {\it nondegenerate} variables and last $(N-s)$ ones {\it singular} variables. Using this geometry we state our main result, an intrinsic form of Harnack inequality.
 \begin{theorem}\label{harnackTHM}
  Let $u$ be a non-negative, bounded, local weak solution to \eqref{E1},\eqref{E1-structure},\eqref{mon-energy}. Let $x_o \in \Omega$ be a point such that $u(x_o)>0$ and $\rho>0$ small enough to allow the inclusion \begin{equation} \label{M} Q_{\mathcal{M},\rho}(x_o) \subseteq \Omega, \quad \quad \text{being} \quad \mathcal{M}= ||u||_{L^{\infty}(\Omega)}^{(2-p)/2}\rho^{\frac{p}{2}}.\end{equation} \noindent
  Assume also that
  \begin{equation}\label{chi}
      \chi:=p+(N-s)(p-2)>0.
  \end{equation}
  
  \noindent Then there exist positive constants $K>1,\bar{\delta}_o\in (0,1)$ depending only on the data such that either
  \begin{equation}\label{alternative-harnack}
      u(x_o) \leq K \rho,
  \end{equation} \noindent or
  \begin{equation}\label{Harnack}
      u(x_o) \leq K \inf_{Q_{\theta,\rho}(x_o)} u\, , \quad \quad  \text{with} \quad \theta= \bar{\delta}_o u(x_o)^{\frac{2-p}{2}} \rho^{\frac{p}{2}}.
  \end{equation}
 \end{theorem} \noindent
%  \begin{remark}
%  {Monotonicity property \eqref{mon-energy} is a natural assumption for existence theory (see for instance \cite{Lions}) and it permits to consider truncations $(u-k)_{\pm}$ for $k \in \R$, as sub-solutions although the equation \eqref{E1} is non homogeneous.}
%  \end{remark}
Condition \eqref{chi} expresses the range of exponents available for the result to hold relatively to the weighted effect of singular and nondegenerate operators into play. When $s$ decreases the range becomes tighter to the parabolic isotropic range for Harnack inequality to hold, i.e. $2N/(N+1)<p<2$. But when $s$ increases the effect of regularisation is stronger and this interval expands until it reaches $1<p<2$.
\vskip0.1cm \noindent The above Theorems \ref{holderTHM} and \ref{harnackTHM} are consequences of the following ones, which are worth of interest on their own. \vskip0.2cm \noindent 
 We prove indeed the following shrinking property, which is typical of both singular parabolic equations (\cite{DB} Lemma 5.1 Chap. IV) and isotropic elliptic equations (\cite{DB2} Prop 5.1 Chap X). 
 \begin{theorem} \label{shrinkingTHM}
  Let $\bar{x} \in \Omega$ and let $u$ be a nonnegative, bounded, local weak solution to \eqref{E1},\eqref{E1-structure},\eqref{mon-energy}. Suppose that for a point $\bar{x} \in \Omega$ and numbers $M,\rho>0$ and $\nu \in (0,1)$ it holds
  \begin{equation}\label{measure-estimate}
      |[u\leq M] \cap Q_{\theta,\rho}(\bar{x}) |\leq (1-\nu) |Q_{\theta,\rho}(\bar{x})|, \quad \text{for} \quad \theta= \rho^{\frac{p}{2}} (\delta M)^{\frac{2-p}{2}},
  \end{equation} \noindent and $Q_{2\theta,2\rho}(\bar{x}) \subset \Omega$, for a number $\delta=\delta(\nu) \in (0,1)$. Then there exist constants $K>1$ and $\delta_o \in (0,1)$ depending only on the data and $\nu$ such that either
  \begin{equation}\label{alternative-shrinking}
    M \leq K \rho,
  \end{equation} \noindent or for almost every  $ x \in Q_{\eta,2\rho}(\bar{x})$ we have  
  \begin{equation}\label{shrinking}
    u(x) \ge \delta_o M/2, \quad \quad \quad  \text{where} \quad \quad \eta= (2\rho)^{\frac{p}{2}} (\delta_o M)^{\frac{2-p}{2}}.   
  \end{equation}
 \end{theorem}\noindent \vskip0.2cm \noindent
 The Theorem above provides also an important expansion of positivity along the singular variables, conceived and modeled in a similar fashion than in (\cite{DGV-Annali}, Theorem 2.3). \vskip0.2cm \noindent Another fundamental tool for our analysis of local regularity is the following integral estimate, which can be seen as an Harnack estimate within the $L^1-L^{\infty}$ topology, and is typical of singular parabolic equations (see for instance \cite{DB}, Prop. 4.1 Chap VII).
 \begin{theorem}\label{l1-linftyTHM}
 Let $u$ be a nonnegative, bounded, local weak solution to \eqref{E1},\eqref{E1-structure},\eqref{mon-energy}. Fix a point $\bar{x} \in \Omega$ and numbers $\theta,\rho>0$ such that $Q_{8\theta,8\rho}(\bar{x}) \subset \Omega$. Then there exists a positive constant $\gamma$ depending only on the data such that either
 \begin{equation} \label{either}
    %  \theta \leq \rho, \quad \text{i.e.} \quad 
     \bigg( \frac{\theta^2}{\rho^p}\bigg)^{\frac{1}{2-p}} \leq \rho,
 \end{equation}\noindent or
 \begin{equation}\label{l1-l1}
      \dashint \dashint_{Q_{\theta,\rho}(\bar{x})} u \, dx \leq \gamma \bigg{\{}\inf_{B_{\frac{\theta}{2}}(\bar{x}')}  \bigg( \dashint_{B_{2\rho}(\bar{x}'')} u(\cdot, x'')\, dx'' \bigg)^{\frac{p}{\chi}}+ \bigg(\frac{\theta^2}{\rho^p} \bigg)^{\frac{1}{2-p}}     \bigg{\}}.
 \end{equation} \noindent If additionally property \eqref{chi} holds, then either we have \eqref{either} or
 \begin{equation} \label{l1-linfty}
     \sup_{Q_{\frac{\theta}{2},\frac{\rho}{2}}(\bar{x})} u \leq \gamma \bigg{\{} \bigg( \frac{\rho^p}{\theta^2} \bigg)^{\frac{N-s}{\chi}} \inf_{B_{\frac{\theta}{2}}(\bar{x}')}   \bigg( \dashint_{B_{2\rho}(\bar{x}'')} u(\cdot, x'') dx''  \bigg)^{\frac{p}{\chi}} +\bigg(  \frac{\theta^2}{\rho^p} \bigg)^{\frac{1}{2-p}}  \bigg{\}}.
 \end{equation}
   \end{theorem}
\noindent

 \subsection{Novelty and Significance}
    Considering fully anisotropic equations as \eqref{prototipo}, a standard statement of regularity requires a bound on the sparseness of the powers $p_i$s. Indeed, in general, weak solutions can be unbounded, as proved in \cite{Giaq}, \cite{Marce0}. We refer to the surveys \cite{Marce4}, \cite{Ming} for an exhaustive treatment of the subject and references. The problem of regularity for anisotropic operators behaving like \eqref{full-EQ} with measurable and bounded coefficients remains a mayor challenge after more than fifty years. Recently some progresses have been made in the parabolic prototype case, as for instance in \cite{Brasco2},\cite{Ural} about Lipschitz continuity, \cite{Ciani} about intrinsic Harnack estimates and \cite{Vazquez} in the singular case for Barenblatt-type solutions. Moreover, as we will see, various parabolic techniques have been applied, but in no circumstance Harnack estimates have been found when more than one spatial dimension was considered. This is due to the fact that usual parabolic techniques rely on the particular structure of a first derivative in time, and are not suitable to manage stronger anisotropies. With the present work, limited to the case $p_i=2$ for $i=1,\dots, s$ and $p_i=p<2$ for $i \in {s+1,\dots,N}$ we are able to prove a purely elliptic pointwise Harnack estimate when the operator acts on the first $s$ variables. Furthermore, we have now a way to understand how these estimates degenerate when $s$ varies; describing, roughly speaking, when the operator is closer to the $p$-Laplacian or to an uniformly elliptic operator (see the discussion after Theorem \ref{harnackTHM}).
    
     \vskip0.2cm \noindent
    In the present work we are interested in bounded solutions, therefore leaving the problem of boundedness to the already rich literature. Our aim is to manage the anisotropic behavior of the operator interpreting its action in correspondence with a suitably adaptd version of the technique developed by E. DiBenedetto (see the original paper \cite{DB-Chen} or the books \cite{DB}, \cite{Urbano}) in order to restore the homogeneity of the parabolic p-Laplacian. Indeed, because of the double derivative, equation \eqref{prototipo} has a wilder heterogeneity of the operator than the parabolic $p$-Laplacian, and 
    % because the parabolic p-Laplacian has only one derivative in time, 
   the intrinsic geometry will be set up according to the order and power of derivatives resulting in the dimensional analysis of the equation.  \vskip0.2cm
   
   \noindent An interesting attempt in this direction has already been done by the some of the authors in \cite{Liao}, in the case of only one nondegenerate variable (see also \cite{LS}). There an expansion of positivity is provided by applying an idea from \cite{DBGV-pams}, shaped on a proper exponential change of variables. Nevertheless, the change of variables in consideration is a purely parabolic tool, so that it does not allow the authors to go through more than one nondegenerate variable. The present work is conceived to fill this gap and to spread new light on the link between classical logarithmic estimates and anisotropic operators.\vskip0.2cm \noindent 
The two fundamental tools that we derive in our work are Theorems \ref{shrinkingTHM}, \ref{l1-linftyTHM}. Theorem \ref{l1-linftyTHM} consists in a $L^{1}-L^{\infty}$ Harnack inequality, which is independent of the other two Theorems, although its proof relies as well on logarithmic estimates \eqref{energy2}. The name $L^{1}-L^{\infty}$-{\it inequality} refers to the fact that it is possible to control the supremum of the function through some $L^1$- integral norm of the function itself. This precious inequality can be used in turn to derive in straight way the H\"older continuity of solutions (see for instance \cite{Io1} for a simple proof in the parabolic setting). \newline
On the other hand Theorem \ref{shrinkingTHM} provides both a shrinking property and an expansion of positivity. The special feature of Theorem \ref{shrinkingTHM} called {\it shrinking} property consists in the fact that from a whatever upper bound to the relative measure of some super-level of the solution it is possible to recover a pointwise estimate of positivity. Its proof is a proper consequence of logarithmic estimates \eqref{energy2} and a suitable choice of test functions (see functions $f$ in Step 1 of the proof of Lemma \ref{SH2-Lemma}). The {\it expansion of positivity} property refers to the possibility to expand along the space (in singular variables) the lower bound yet gained. The proof of this property is therefore linked to the measure theoretical approach of this shrinking property, and it is an adaptation of an idea of \cite{DGV-Annali}. Here to end the proof of Theorem \ref{shrinkingTHM} we use in a crucial way the shrinking property to reach a critical mass and use Lemma \ref{DG}.\newline
Finally, in order to prove Theorem \ref{harnackTHM} we use an argument originally conceived by Krylov and Safonov in \cite{KS} to reach a certain controlled bound on the solution in terms of the solution itself, and then use repeatedly Theorem \ref{l1-linftyTHM} to achieve an upper bound on the measure of some super-level set of the solution. Nonetheless, the argument of Krylov and Safonov gave us this information around an unknown point. Therefore we apply Theorem \ref{shrinkingTHM} to expand the positivity until the desired neighborhood of the initial point and get the job done. The idea is an adaptation of the techniques originally developed in \cite{DGV-Annali} to the case of anisotropic elliptic equations \eqref{E1}-\eqref{E1-structure}.

    \subsection{Structure of the paper} \noindent 
    In Section \ref{preliminaries} we recall major functional tools and use them to derive fundamental properties of solutions as energy estimates, logarithmic estimates and some integral estimates. Then in Section \ref{3} we prove Theorem \ref{shrinkingTHM}, in Section \ref{4} we prove H\"older continuity of solutions while in Section \ref{5} we prove the $L^1-L^{\infty}$ Harnack estimate \eqref{l1-linfty}. Section \ref{7} is devoted to the proof of Theorem \ref{harnackTHM}. Technicalities and standard material has been collected in a final section, Section \ref{8}, to leave space along the previous text to what is really new.\vskip0.2cm 
    \noindent {\bf Notations}: 
\begin{itemize}

\item[-]
If $\Omega$ is a measurable subset of $\R^{N}$, we denote by $|\Omega|$ its Lebesgue measure. We will write $\Omega \subset \subset \R^N$ when $\Omega$ is an open bounded set.
\item[-]
For $r>0$, $\bar{x}=(\bar{x}',\bar{x}'') \in \R^{s} \times \R^{N-s}$, we denote by $B_{r}(\bar x)$ the ball of radius $r$ and center $\bar{x}$; the standard polydisc is denoted by $Q_{\theta,\rho}=B_{\theta}(\bar{x}') \times B_{\rho}(\bar{x}'') \subset \R^N$. Furthermore, by $w_s=|B_1(0')|$ and $w_{N-s}=|B_1(0'')|$ we denote the measures of the respective unit balls.
\item[-]
The symbol $\forall_{\text{ae}}$ stands for {\it -for almost every-} .   
\item[-]
For a measurable function $u$, by $\inf u$ and $\sup u$ we understand the essential infimum and supremum, respectively;  when $u:\Omega \to \R$  and $a\in \R$, we  omit the domain when considering sub/super level sets, letting $\big[u\gtreqless a\big]=\big\{x\in E: u(x)\gtreqless a\big\}$; if $u$ is defined on some open set $\Omega\subset \R^N$, we let  $\partial_iu=\frac{\partial }{\partial x_i}u$ denote the distributional derivatives.
\item[-] For numbers $B,C>0$ we write $C\wedge B= \max \{ B,C\}$.
\item[-] We make the usual convention that a constant $\gamma>0$ depending only on the data, i.e. $\gamma=\gamma(N, 2, p, C_1, C_2, C)$, may vary from line to line along calculations.
\end{itemize}

 {\centering{\section{Preliminaries} \label{preliminaries}}} \noindent In this Section we collect the basic tools that will be used along the overall theory. For the sake of readability, simpler and well-known proofs are postponed to the Appendix (Section \ref{8}), while most relevant passages that bring to light our method are detailed and highlighted. 

\vskip0.2cm \noindent
\subsection{Functional and standard tools}
We recall the embedding $W_o^{1,\bf{p}}(\Omega) \hookrightarrow W^{1,\bar{p}^*}_o(\Omega)$ proved by M. Troisi in \cite{Troisi}.
\begin{lemma}  \label{Troisi}
Let $\Omega \subset \subset \R^N$ and consider a function $u \in W_o^{1,\bf{p}}(\Omega)$, $p_i>1$ for each $i \in \{1,..,N\}$. Assume $\bar{p}<N$ and let 
\begin{equation}\label{sobolevp}
\bar{p}^{*}= \frac{N \bar{p}}{N-\bar{p}}.
\end{equation} \noindent Then there exists a positive constant $\gamma(N,\bar{p})$ such that
\begin{equation}\label{troisi}
    ||u||_{L^{\bar{p}^*}(\Omega)}^N \leq \gamma \prod_{i=1}^N ||\partial_i u ||_{L^{p_i}(\Omega)}.
\end{equation}
\end{lemma} \noindent It is worth pointing out that without vanishing initial datum this embedding fails in general (see \cite{KK},\cite{Schmeiser} for counter-examples). A simple calculation reveals that condition $\bar{p} <N$ is always in force in our case in study. Next Lemma introduces a well-known weighted Poincaré inequality (see for instance Prop.2.1 in \cite{DB}), that will be useful when estimating the logarithmic function.

\begin{lemma}
Let $B_{\rho}$ be a ball of radius $\rho>0$ about the origin, and let $\varphi \in C(B_{\rho})$ satisfy $0 \leq \varphi(x) \leq 1$ for each $x \in B_{\rho}$ together with the condition that the level sets $[\varphi >k] \cap B_{\rho}$ are convex for each $k \in (0,1)$. Let $g \in W^{1,p}(B_{\rho})$, and assume that the set 
\[
\mathcal{E} = [g=0] \cap [\varphi=1]
\] has positive measure. The there exists a constant $\gamma>0$ depending only upon $N,p$ such that 
\begin{equation}\label{WP}
    \int_{B_{\rho}} \varphi |g|^p\, dx\leq \gamma  \rho^p \bigg( \frac{|B_{\rho}|}{|\mathcal{E}|} \bigg)^p \int_{B_{\rho}} \varphi |D g |^p \, dx.
\end{equation}
\end{lemma}
\noindent

\subsection{Properties of solutions to \eqref{E1}-\eqref{E1-structure}}
The following classical {\it Energy Estimates} can be proved by a standard choice of test functions.\vskip0.1cm \noindent 
\begin{lemma}
Let $u$ be a bounded local weak solution to \eqref{E1} with structure conditions \eqref{E1-structure}. Then there exists a positive constant $\gamma$ such that for any polydiscs $Q_{\theta, \rho}(\bar{x}) \subseteq \Omega$, any $k \in \R$ and any $\zeta \in C_o^{\infty}(Q_{\theta, \rho}(\bar{x}))$ such that $0 \leq \zeta \leq 1$ it holds
\begin{equation}\label{energy1}
\begin{aligned}
&\sum_{i=1}^s \int \int_{Q_{\theta, \rho}(\bar{x})} | \partial_i (u-k)_{\pm}|^{2} \zeta^2 \, dx + \sum_{i=s+1}^N\int \int_{Q_{\theta, \rho}(\bar{x})} | \partial_i (u-k)_{\pm}|^{p} \zeta^2 \, dx\leq\\
&  \gamma  \int \int_{Q_{\theta, \rho}(\bar{x})} \bigg{\{}  |(u-k)_{\pm}|^{2}  \sum_{i=1}^s |\partial_i \zeta|^2   + | (u-k)_{\pm}|^{p} \sum_{i=s+1}^N  |\partial_i \zeta|^p  + C^p \chi_{[(u-k)_{\pm}>0]}        \bigg{\}}\, dx .
\end{aligned}
\end{equation}
\end{lemma} \noindent See Section \ref{8} for the classical proof. Next Lemma is a sort of measure theoretical maximum principle. It asserts that if a certain sub-level set of the solution reaches a critical mass, then the solution is above a multiple of the level on half sub-level set. We agree to refer to it as usual in literature by the epithet {\it Critical Mass} Lemma ({\it De Giorgi-type} Lemma is used equivalently). We state it just for sub-level sets, a similar statement being true for super-level sets.

\begin{lemma} \label{DG}
Let $\bar{x} \in \Omega$ and $\theta,\rho>0$ such that $Q_{4\theta,4\rho}(\bar{x}) \subseteq \Omega$. Let $\mu^+,\mu^-,\omega$ be nonnegative numbers such that
\[
\mu^+ \ge \sup_{Q_{\theta,\rho}(\bar{x})} u, \quad \quad \mu^- \leq \inf_{Q_{\theta,\rho}(\bar{x})} u, \quad \quad  \omega \ge \mu^+-\mu^-.
\]
Now, let $u$ be a bounded function satisfying the energy estimates \eqref{energy1} and fix $a, \xi \in (0,1)$. Then there exists a number $\nu\in (0,1)$ whose dependence from the data is specified by \eqref{nu} and such that if
\begin{equation}\label{DG-HP}
    |[u\leq \mu^-+\xi \omega] \cap Q_{\theta,\rho}(\bar{x})| \leq \nu | Q_{\theta,\rho}|,
\end{equation} \noindent then either $\xi \omega \leq \rho$ or

\begin{equation}\label{DG-TH}
    u \ge \mu^-+a\xi \omega, \quad \quad \forall_{\text{ae}} \, \,  x \in Q_{\theta/2, \rho/2}(\bar{x}).
\end{equation}

\end{lemma}

\begin{proof}
We suppose without loss of generality that $\bar{x}=0$. For $j=0,1,2..$ let us set
\begin{equation}
    \begin{cases}
    k_j= \mu^- + a \xi \omega + \frac{(1-a)\xi \omega}{2^j},\\
    \rho_j= \frac{\rho}{2}+\frac{\rho}{2^{j+1}}, \quad \theta_j= \frac{\theta}{2}+ \frac{\theta}{2^{j+1}},
    \end{cases} \quad A_j= Q_{\theta_j,\rho_j} \cap [u <k_j],
\end{equation}\noindent and let $\zeta_j \in C_o^{\infty}(Q_j)$ be a cut-off function between $Q_j$ and $Q_{j+1}$ such that $\zeta_j \equiv 1$ on $Q_{j+1}$,$0 \leq \zeta_j \leq 1$ and therefore satisfying
\[
|\partial_i \zeta_j| \leq \frac{2^{j+2}}{\theta} \quad \forall i \in \{1,..,s\}, \quad \text{\&} \quad  |\partial_i \zeta_j| \leq \frac{2^{j+2}}{\rho} \quad \forall i \in \{s+1,..,N\}.
\]
Then, combining a precise use of H\"older inequality and Troisi's embedding \eqref{troisi} to the energy estimates \eqref{energy1}, leads us to the estimate
\begin{equation}\label{DGa}
    \begin{aligned}
    \bigg(\frac{(1-a) \xi \omega}{2^{j+1}}  \bigg)^{\bar{p}}& |A_{j+1}| \leq \int \int_{Q_j} (u-k)_{-}^{\bar{p}} \zeta_j^{\bar{p}} \, dx \leq\\
    &\bigg( \int \int_{Q_j} [ (u-k)_{-} \zeta_j]^{\frac{N \bar{p}}{N-\bar{p}}} \bigg)^{\frac{N-\bar{p}}{N}} |A_j|^{\frac{\bar{p}}{N}} \leq\\
    & \bigg[ \prod_{i=1}^N \bigg( \int \int_{Q_j} |\partial_i (u-k)_{-}|^{p_i} \, dx \bigg)^{\frac{1}{p_i}} \bigg]^{\frac{\bar{p}}{N}}|A_j|^{\frac{\bar{p}}{N}} \leq \\
    & \bigg[ \int \int_{A_j} \bigg{\{}\sum_{i=1}^s | \partial_i (u-k)_{\pm}|^{2} \zeta^2 + \sum_{i=s+1}^N | \partial_i (u-k)_{\pm}|^{p} \zeta^2 + C^p\bigg{\}} dx \bigg]|A_j|^{\frac{\bar{p}}{N}} \leq \\
    % & \gamma \bigg{\{}  \xi^2 \omega^2 \bigg[a+\bigg( \frac{1-a}{2^j}\bigg)  \bigg]^2 \frac{2^{2j}}{\theta^2} + \xi^p \omega^p \bigg[a+\bigg( \frac{1-a}{2^j}\bigg)  \bigg]^p  \frac{2^{pj}}{\rho^p}\bigg{\}}|A_j|^{1+\frac{\bar{p}}{N}} \leq \\
    &\gamma 2^{2j} \bigg{\{}\frac{(\xi \omega)^p}{\rho^p} \bigg[ 1+\frac{(\xi \omega)^{2-p} \rho^p}{\theta^2}+ \bigg(\frac{C \rho }{\xi \omega}\bigg)^p \bigg] \bigg{\}} |A_j|^{1+\frac{\bar{p}}{N}}.
    \end{aligned}
\end{equation} \noindent By assumption $\xi \omega >\rho$ the third term on right hand side is smaller than $1$. If we define $Y_j= |A_j|/ |Q_j|$, we divide \eqref{DGa} by $|Q_{j+1}|$ and we observe that $|Q_{j}| \leq \gamma 2^j |Q_{j+1}|\approx (\theta^s \rho^{N-s})$, previous estimate can be written as
\begin{equation*}
\begin{aligned}
    Y_{j+1} \leq& 
    \quad \gamma \, \frac{2^{(2+p)j} (1-a)^{-\bar{p}} (\xi \omega)^{p-\bar{p}}}{\rho^p}\bigg[ 1+\frac{(\xi \omega)^{2-p} \rho^p}{\theta^2} \bigg] (\theta^s \rho^{N-s})^{\frac{\bar{p}}{N}}\,  Y_j^{1+\frac{\bar{p}}{N}}  \\
    & \leq \gamma 2^{(2+p)j} (1-a)^{-\bar{p}} \bigg( \frac{\theta}{\rho^{p/2} (\xi \omega)^{\frac{2-p}{2}}} \bigg)^{\frac{s \bar{p}}{N}} \bigg[ 1+\frac{(\xi \omega)^{2-p} \rho^p}{\theta^2} \bigg]   Y_j^{1+\frac{\bar{p}}{N}},
    \end{aligned}
\end{equation*} \noindent
by simple manipulation on the various exponents. We evoke Lemma \ref{Fast} to declare that if
\begin{equation}\label{nu}
\begin{aligned}
    Y_o= \frac{|[u \leq \mu^-+\xi  \omega]  \cap Q_{\theta,\rho}|}{|Q_{\theta,\rho}|} \leq 
    \gamma^{\frac{-N}{\bar{p}}}(1-a)^N \bigg( \frac{\theta}{\rho^{p/2} (\xi \omega)^{\frac{2-p}{2}}}  \bigg)^{-s} \bigg[ 1+\frac{(\xi \omega)^{2-p} \rho^p}{\theta^2} \bigg]^{-\frac{N}{\bar{p}}}=: \nu\, ,
    \end{aligned}
\end{equation}\noindent then $Y_j \rightarrow 0$ for $j \rightarrow \infty$ and the proof is concluded by specifying $Y_{\infty}=\lim_{j \rightarrow \infty} Y_j=0$.
\end{proof} 
\begin{remark}\label{R-geometry}
We observe that within geometry \eqref{cylinders} the choice $ \theta= \rho^{\frac{p}{2}} (\xi \omega)^{\frac{2-p}{2}}$ sets $\nu$ free from any other dependence than the initial data.
\end{remark} \noindent The following Lemma estimates the essential supremum of solutions by quantitative integral averages of the solution itself. Its proof is similar to the one of (\cite{Liao}, Prop. 8) and it is postponed to the Appendix.

\begin{lemma} \label{DG-IT}
Le $u$ be a locally bounded local weak solution to \eqref{E1}-\eqref{E1-structure}. Let $ 1\leq l \leq 2$ and 
\begin{equation}\label{lambdas}
    N(\bar{p}-2)+l\bar{p}>0.
\end{equation} \noindent 
Then there exist constants $\gamma,C>1$ depending only on the data, such that for all polydiscs $Q_{2\theta,2\rho}\subset \Omega$ we have either
\begin{equation} \label{raggae0}
\bigg( \frac{\theta^2}{\rho^p} \bigg)^{\frac{1}{2-p}}\leq C \rho, \quad \text{or}
\end{equation}
\begin{equation}\label{Ls-Linfty}
    \sup_{Q_{\theta/2},\rho/2} u  \leq \gamma  \bigg( \frac{\rho^p}{\theta^2} \bigg)^{( \frac{N-s}{p}) (\frac{\bar{p}}{N(\bar{p}-2)+l\bar{p}})}  \bigg( \dashint  \dashint_{Q_{\theta,\rho}} u_+^l\, dx  \bigg)^{\frac{\bar{p}}{N(\bar{p}-2)+l\bar{p}}} +\gamma \bigg(\frac{\theta^2}{\rho^p} \bigg)^{\frac{1}{2-p}}.
\end{equation}\noindent 
\end{lemma}

\noindent Note that \eqref{lambdas} with $l=1$ corresponds to \eqref{chi}. Finally we give detailed description of the main analytical tool of the present work, the following {\it Logarithmic Estimates}.

\begin{lemma} \label{log-estimates}
Let $u$ be a bounded weak solution to \eqref{E1},\eqref{E1-structure},\eqref{mon-energy}. Then for any $Q_{2\theta,2\rho}(\bar{x}) \subset \Omega$, any $k \in \R$ and any function $f \in C^1(\R;\R_+)$ with $f'>0$, there exists a constant $\gamma>1$ such that the following estimate holds for each $0<t<\theta$, and each $\zeta \in W^{1,p}_o(B_{2\rho}(\bar{x}''))$,

\begin{equation}\label{energy2}
    \begin{aligned}
     &\int \int_{Q_{t,\rho}(\bar{x})} f'((u-k)_{\pm}) (t^2-|x'-\bar{x}'|^2) \zeta^p(x'') \bigg{\{} \sum_{i=1}^s|\partial_i (u-k)_{\pm}|^2+ \sum_{i=s+1}^N |\partial_i (u-k)_{\pm}|^p   \bigg{\}}dx\\
    &  \leq \gamma t^{s+1} \frac{d}{\d t}\bigg( |B_{t}(\bar{x}')|^{-1} \int \int_{Q_{t,\rho}(\bar{x})} \bigg[ \int_0^{(u-k)_{\pm}} f(\tau) d\tau    \bigg] \zeta^p(x'') dx \bigg)+\\
     &+ \gamma \int \int_{Q_{t,\rho}(\bar{x})}  (t^2-|x'-\bar{x}'|^2) \bigg{\{} \sum_{i=s+1}^N
     \frac{[f((u-k)_{\pm})]^p}{[f'((u-k)_{\pm})]^{p-1}}  ( | \partial_i \zeta|^p+\zeta^p) + f'((u-k)_{\pm}) \bigg{\}}  dx.
    \end{aligned}
\end{equation}
\end{lemma}

\begin{proof}
We test the equation \eqref{mon-energy} with a nonnegative function
\[\psi(x)= f(\pm(u(x)-k)_{\pm}) \phi(x') \zeta^p(x'') \in W_o^{1,[2,p]}(Q_{\theta,\rho}),\] 
being $\phi \in W^{1,2}_o(B_{2\theta}(\bar{x}'))$ a test function, and we use the structure conditions \eqref{E1-structure} to get 
% the following formal algebraic computations under the sign of integral: first the test gives
% \begin{align*}
%   0= & \int \int \bigg{\{} \sum_{i=1}^s (\partial_i u) \bigg[ f'((u-k)_{\pm}) (\partial_i (u-k)_{\pm}) \phi(x') \zeta^p+ f((u-k)_{\pm}) (\partial_i \phi ) \zeta^p \bigg]+\\
%     & +\sum_{i=s+1}^N A_i(x, \nabla u) \bigg[ f'((u-k)_{\pm}) (\partial_i (u-k)_{\pm}) \phi \zeta^p+ p f((u-k)_{\pm}) \phi (\partial_i \zeta ) \zeta^{p-1} \bigg]\bigg{\}} dx,
% \end{align*} \noindent  Using structure conditions \eqref{E1-structure} we get
\begin{align*}
    &\int \int  f'((u-k)_{\pm}) \phi \zeta^p \bigg{\{}\sum_{i=1}^s |\partial_i (u-k)_{\pm}|^2   + C_1 \sum_{i=s+1}^N |\partial_i (u-k)_{\pm}|^p \bigg{\}}\, dx \leq \\
    & \int \int \bigg{\{}  
    \sum_{i=s+1}^N C f'((u-k)_{\pm}) |\partial_i (u-k_{\pm}) |\,  \phi \zeta^p +\sum_{i=1}^s (- \partial_i \phi) f((u-k)_{\pm}) (\partial_i (u-k)_{\pm}) \zeta^p+\\
     &\quad \quad \quad\quad \quad \quad\quad 
     \quad \quad \quad  \quad  \quad \quad 
    + \bigg[ C_2 |\partial_i (u-k)_{\pm}|^{p-1}+C \bigg] f((u-k)_{\pm}) \phi \zeta^{p-1} |\partial_i \zeta| \bigg{\}} \, dx,
\end{align*} \noindent 
the integrals being taken over the polydisc $Q_{\theta,\rho}(\bar{x})$. We use repeatedly Young's inequality to get, for the first term on the right
\[
Cf' |\partial_i (u-k_{\pm}) |\,  \phi \zeta^p \leq  f' \phi \bigg( \epsilon  |\partial_i (u-k)_{\pm} |^p \zeta^p + \gamma(\epsilon) C^{\frac{p}{p-1}} \zeta^p   \bigg),
\]while the third term on the right is estimated with
\begin{align*}
|\partial_i (u-k)_{\pm}|^{p-1} f  \phi \zeta^{p-1} \partial_i \zeta (f')^{\frac{p-1}{p}}(f')^{\frac{1-p}{p}}\leq \phi \bigg(\epsilon f' |\partial_i (u-k)_{\pm}|^p \zeta^p + C(\epsilon)  | \partial \zeta|^p f^p(f')^{1-p}    \bigg).
\end{align*} \noindent Similarly fourth term is estimated by
\[
C f \phi \zeta^{p-1} |\partial_i \zeta| (f')^{\frac{p-1}{p}}(f')^{\frac{1-p}{p}}\, \leq  \phi \bigg( \epsilon  f^p(f')^{1-p} |\partial_i \zeta |^p +\gamma(\epsilon)  C^{\frac{p}{p-1}} f' \zeta^p \bigg).
\] 
 \noindent Gathering all the pieces together and choosing accordingly $\epsilon$ small enough to reabsorb on the left the energy terms, we get 
\begin{align*}
    \int \int&  f'((u-k)_{\pm}) \phi \zeta^p \bigg[\sum_{i=1}^s |\partial_i (u-k)_{\pm}|^2   +  \sum_{i=s+1}^N |\partial_i (u-k)_{\pm}|^p \bigg]\, dx \leq\\
    &\gamma \int \int \bigg[ 
    f'((u-k)_{\pm}) \phi \zeta^p- \sum_{i=1}^s (\partial_i \phi) (\partial_i (u-k)_{\pm}) f((u-k)_{\pm}) \zeta^p \bigg] dx+\\
    &+ \gamma \int \int \bigg[\sum_{i=s+1}^N \frac{f^p((u-k)_{\pm})}{(f'(u-k)_{\pm})^{p-1}} \phi ( |\partial_i \zeta|^p+|\zeta^p|) + f'((u-k)_{\pm}) \phi \zeta^p \bigg] \,dx.
\end{align*} \noindent Finally, choose $0 \leq \zeta(x'') \leq 1$ and $\phi(x')=(t^2-|x'-\bar{x}'|^2)_{+}$ and estimate the second term on the right by use of Green-Ostrogradsky's formula with

\begin{equation} \label{Stokes}
\begin{aligned}
\int &\int_{Q_{t,\rho}(\bar{x})} f((u-k)_{\pm}) \zeta^p(x'')  \sum_{i=1}^s (\partial_i (u-k)_{\pm}) 2(x_i'-\bar{x}_i') \, dx =\\
& 2 \int \int_{Q_{t,\rho}(\bar{x})} \sum_{i=1}^s\bigg[ \partial_i \bigg( \int_0^{(u-k)_{\pm}} f(\tau) d\tau \cdot  (x_i'-\bar{x}_i')  \bigg) -2s \bigg( \int_0^{(u-k)_{\pm}} f(\tau) d \tau \bigg)\bigg]  \zeta^p dx =\\
&2t \int_{\partial B_t(\bar{x}')} \int_{B_{\rho}(\bar{x}'')} \bigg( \int_0^{(u-k)_{\pm}} f(\tau) d \tau \bigg)\zeta^p\, dx-2s \int \int_{Q_{t,\rho}(\bar{x})} \bigg( \int_0^{(u-k)_{\pm}} f(\tau) d \tau \bigg)\bigg]  \zeta^p dx \\
&= 2 t^{s+1} \frac{d}{dt} \bigg[t^{-s} \int \int_{Q_{t,\rho}(\bar{x})} \bigg( \int_0^{(u-k)_{\pm}} f(\tau) d \tau \bigg)\zeta^p(x'') dx   \bigg],
\end{aligned}
\end{equation} \noindent and the proof is concluded.
\end{proof} \noindent

\vskip0.2cm

{\centering{
\section{Proof of Theorem \ref{shrinkingTHM}} \label{3}}}
\noindent 
We start by proving two main Lemmas, whose combination will provide an easy proof of Theorem \ref{shrinkingTHM}. The first one turns a measure estimate given on the intrinsic polydisc into a measure estimate on each $(N-s)$-dimensional slice.
\vskip0.2cm \noindent
\begin{lemma} \label{SH1-Lemma}
Let $u$ be a bounded weak solution to \eqref{E1}, with structure conditions \eqref{E1-structure}, \eqref{mon-energy}. Assume that for some $M,\rho>0$ and some $\alpha\in (0,1)$ the following estimate holds
\begin{equation}\label{S1-HP}
    |[u(\cdot) \leq M] \cap Q_{\theta,\rho}(\bar{x})| \leq (1-\alpha) |Q_{\theta,\rho}|, \quad \quad \theta= \rho^{\frac{p}{2}} (\delta M)^{\frac{2-p}{2}},
\end{equation}\noindent for some $\delta < \alpha^4(1+\alpha)^s / (1+\alpha^2)$. Then there exist numbers $s_1>1$, $\delta_1 \in (0,1)$ depending only on the data and $\alpha$ such that
\begin{equation}\label{alternative}
    M\leq \rho , \quad \quad \text{or}
\end{equation}\noindent
\begin{equation}\label{S1-TH}
    |[u(y',\,\cdot) \leq 2^{-s_1}M]\cap B_{\rho}(\bar{x}'')| \leq (1-\alpha^4/2) |B_{\rho}(\bar{x}'')|, \quad \quad \forall_{\text{ae}}\,\, y' \in B_{(1-\delta_1)\theta}(\bar{x}').
\end{equation}
\end{lemma}
\vskip0.1cm

\begin{proof}
Let $\sigma, \delta_1 \in (0,1)$ to be chosen later, let $\zeta(x'') \in C_o^{\infty}(B_{\rho}(\bar{x}''))$ be such that 
\begin{equation*} \begin{cases} 0 \leq \zeta \leq 1,\\ \zeta(x'')\equiv 1, &  x''\in B_{(1-\sigma) \rho}(\bar{x}''), \end{cases} \quad \quad \text{\&}  \quad \quad |\partial_i \zeta| \leq \frac{\gamma}{\sigma \rho}, \quad \forall i=s+1,..,N.\end{equation*} \noindent 
Let $y' \in B_{(1-\delta_1) \theta}(\bar{x}')$ be a Lebesgue point for the function 
\[
\int_{B_{\rho}(\bar{x}'')} (u(\cdot,x)-M)_{-}^2 \zeta^p(x'')\, dx,
\] let us call $z=(y',\bar{x}'')$ and use just the right hand side of inequality \eqref{energy2},
% \begin{equation*}
%     \begin{aligned}
%     &0\leq t^{s+1}\frac{d}{dt} \bigg[ t^{-s} \int \int_{Q_{t,\rho}(z)} \bigg( \int_0^{(u-k)_{-}} f(\tau) \, d\tau\bigg) \zeta^p(x'') dx \bigg]+\\
%     & + \gamma   \int \int_{Q_{t,\rho}(z)}(t^2-|x'-\bar{x}'|^2) \bigg{\{} \sum_{i=s+1}^N \frac{f^p((u-k)_{-})}{[f'((u-k)_{-})]^{p-1}}  ( | \partial_i \zeta|^p+1) +  f'(u-k)_{-}   \bigg{\}} dx 
%     \end{aligned}
% \end{equation*} \noindent
with $f(u)=u$, $f'\equiv 1$, $k=M$, to get
% \footnote{ According with $(\sigma \rho)^{-1} \ge 1$, and $(t^2-|x'-\bar{x}'|^2)\leq t^2$ .}
\begin{equation*}
    \begin{aligned}
    0\leq \frac{d}{dt} \bigg[ t^{-s} \int& \int_{Q_{t,\rho}(z)} \frac{(u-M)^2_{-}}{2} \zeta^p(x'') dx \bigg]+
    \\
    & + \frac{\gamma t^{-s-1}}{(\sigma \rho)^p} \int \int_{Q_{t,\rho}(z)}t^2 \bigg{\{} \sum_{i=s+1}^N (u-M)_{-}^p  +(\sigma \rho)^p  \bigg{\}} dx.
    \end{aligned}
\end{equation*} \noindent Now we integrate this inequality on $t \in (0, \delta_1\theta)$ and estimate the various terms. The first term can be evaluated by
\begin{equation*}
    \begin{aligned}
   \bigg[ t^{-s}  \int \int_{Q_{t,\rho}(z)}& \frac{(u-M)_{-}^2}{2}\zeta^p(x'') \, dx \bigg]_0^{\delta_1\theta}=\\
    % &\frac{1}{(\delta_1 \theta)^s} \int \int_{Q_{t,\rho}(z)} \frac{(u-N)_{-}^2}{2} \, dx- \lim_{t \rightarrow 0} t^{-s} \int \int_{Q_{t,\rho}(y',\bar{x}'')} \frac{(u-N)_{-}^2}{2} \, dx=\\
    & \frac{1}{(\delta_1 \theta)^s} \int \int_{Q_{\delta_1\theta,\rho}(z)} \frac{(u-M)_{-}^2}{2}\zeta^p(x'') \, dx-  \int_{B_{\rho}(\bar{x}'')} \frac{(u(y',\cdot)-M)_{-}^2}{2}\zeta^p(x'') \, dx,
    \end{aligned}
\end{equation*}

\noindent using that $y'$ is a Lebesgue point. Second term is estimated with 
\begin{equation*}
    \begin{aligned}
&\gamma \int_0^{\delta_1 \theta}\bigg(\frac{t^{-s+1}}{(\sigma \rho)^p)} \int \int_{Q_{t,\rho}(z)} (u-M)_{-}^p  dx\bigg)dt                     \leq 
    % &\gamma \int_0^{\delta_1 \theta} \frac{t^{-s+1}}{(\sigma \rho)^p)} N^p |B_t(y')|\, |B_{\rho}(\bar{x}'')|\, dt \leq\\
    % \frac{\gamma N^p}{(\sigma \rho)^p} |B_{\rho}(\bar{x}'')| \, \int_0^{\delta_1 \theta} t \, dt =
    \frac{\gamma M^p}{(\sigma \rho)^p} (\delta_1 \theta)^2|B_{\rho}(\bar{x}'')|,
    % \\
    % &  \gamma \int_0^{\delta_1 \theta} \bigg( t^{-s+1} \int \int_{Q_{t,\rho}(z)} dx \bigg) dt \leq
    % % &\gamma \int_0^{\delta_1\theta} \bigg( t^{-s+1} |B_t(y')| \bigg)\, |B_{\rho}(\bar{x}'')| \leq\\
    % \gamma (\delta_1 \theta)^2 |B_{\rho}(\bar{x}'')|.
    \end{aligned}
\end{equation*} \noindent and third term similarly. Gathering all together we obtain
\begin{equation*}
    \begin{aligned}
\int_{B_{(1-\sigma)\rho}(\bar{x}'')}& (M-u(y',x''))_{+}^2\, dx \leq\\
& \frac{1}{(\delta_1\theta)^s}  \int \int_{Q_{\delta_1\theta,\rho}(z)} (M-u)_{+}^2\, dx + \gamma (\delta_1 \theta )^2 |B_{\rho}(\bar{x}'')| \bigg{\{} \bigg(\frac{M}{\sigma \rho}   \bigg)^p +1 \bigg{\}}.
\end{aligned} \end{equation*} \noindent 
We observe that $B_{\delta_1}(y') \subset B_{\theta}(\bar{x}')$ and that by imposing the natural intrinsic geometry the term in parenthesis can be ruled. Indeed we let $\theta= \rho^{p/2} (\delta N)^{(2-p)/2}$ and compute
\begin{equation*}
    \begin{aligned}
    M^2\bigg( 1-\frac{1}{2^{s_1}} \bigg)^2& | [u(y',\cdot) \leq \frac{M}{2^{s_1}}] \cap B_{(1-\sigma)\rho}(\bar{x}'')|\leq\\
  & \int_{B_{(1-\sigma)\rho}(\bar{x}'')} (M-u(y',x''))_{+}^2\, dx \leq\\
    & \frac{M^2}{(\delta_1 \theta)^s} |[u \leq M] \cap Q_{\theta,\rho}(\bar{x})| + \gamma M^2 |B_{\rho}(\bar{x}'')| \bigg{\{} \frac{(\delta_1 \theta)^2M^{p-2}}{(\sigma \rho)^p}   + \frac{(\delta_1 \theta)^2}{M^2} \bigg{\}}\leq\\
    & M^2 |B_{\rho}(\bar{x}'')|\,\bigg{\{}  \frac{1}{(\delta_1 )^s} (1-\alpha) +\gamma \bigg(\frac{ \delta^{2-p} \delta_1}{\sigma^p}\bigg)\bigg{\}}, \end{aligned}
\end{equation*}\noindent  using $M \ge \rho$ and the hypothesis \eqref{S1-HP}. We estimate from below in the whole $B_{\rho}(\bar{x}'')$ by 
\begin{equation*}
    \begin{aligned}
    &|[u(y',\cdot) \leq M2^{-s_1}] \cap B_{\rho}(\bar{x}'')| \leq
    % \\
    % & |[u(y',\cdot) \leq N2^{-s_1}] \cap B_{(1-\sigma)\rho}(\bar{x}'')| + | B_{\rho} / B_{(1-\sigma)\rho}| \leq\\
    % &
    |[u(y',\cdot) \leq M2^{-s_1}] \cap B_{(1-\sigma)\rho}(\bar{x}'')|+ (N-s) \sigma |B_{\rho}(\bar{x}'')|.
     \end{aligned}
\end{equation*}\noindent 
% Indeed,
% \[
%     | B_{\rho}/ B_{(1-\sigma)\rho}| = w_{N-s} \{ \rho^{N-s}-[(1-\sigma)\rho]^{N-s}\}= w_{N-s} \rho^{N-s} \{1-(1-\sigma)^{N-s}  \} \leq |B_{\rho}| (N-s) \sigma.
% \]
Combining this remark with previous calculations we have the inequality
\begin{equation*}
    \begin{aligned}
    |[u(y',\cdot)& \leq M2^{-s_1}] \cap B_{\rho}(\bar{x}'')|\leq\\
    % & |[u(y',\cdot) \leq N2^{-s_1}] \cap B_{(1-\sigma)\rho}(\bar{x}'')|+ (N-s) \sigma |B_{\rho}(\bar{x}'')|\leq\\
    & (1-2^{-s_1})^{-2} |B_{\rho}(\bar{x}'')|\,\bigg{\{}  \frac{1}{(\delta_1 )^s} (1-\alpha) + \gamma \bigg(\frac{ \delta^{2-p} \delta_1}{\sigma^p}\bigg)\bigg{\}} + (N-s) \sigma |B_{\rho}(\bar{x}'')|=\\
    &|B_{\rho}(\bar{x}'')| \bigg{\{}\frac{(1-\alpha)}{\delta_1^s (1-2^{-s_1})^2}+\frac{\gamma \delta^{2-p}\delta_1}{\sigma^p(1-2^{-s_1})^2} +(N-s) \sigma \bigg{\}}.
     \end{aligned}
\end{equation*}\noindent To conclude the proof we choose $\delta_1$,$s_1$,$\sigma$,$\delta$, from the conditions
\[ (1-2^{-s_1})^2=1+\alpha^2, \quad  \delta_1^{-s_1}=1+\alpha, \quad (N-s) \sigma= \alpha^4/4, \quad \delta^{2-p} =\frac{\alpha^4}{4}\bigg(\frac{\gamma \delta_1}{\sigma^p(1-2^{-s_1})^2} \bigg)^{-1}. \]  
% We remark that these parameters are chosen only in terms of the data. In conclusion the measure estimate
% \[
% |[u(y',\cdot) \leq N2^{-s_1}] \cap B_{\rho}(\bar{x}'')|\leq|B_{\rho}(\bar{x}'')| \bigg{\{}(1-\alpha)(1+\alpha)(1+\alpha^2)+\frac{\alpha^4}{4}+\frac{\alpha^4}{4} \bigg{\}}= (1-\alpha^4/2) |B_{\rho}(\bar{x}'')|
% \] is valid for each $y' \in B_{(1-\delta_1)\theta}(\bar{x}')$.
\end{proof} \noindent 

\noindent 
Second and following Lemma is what is called in literature a {\it shrinking} Lemma. Indeed, from a given relative measure information on a level set, it allows to shrink as much as we need the relative measure on a lower-level set. Even more interesting, it provides also an expansion of positivity along singular variables.

\begin{lemma} \label{SH2-Lemma}
Let $\bar{x} \in \Omega$ and $u$ be a bounded weak solution to \eqref{E1},\eqref{E1-structure},\eqref{mon-energy}. Let us suppose that for some $M,\rho>0$ and some $\beta \in (0,1)$ the following estimate holds for almost every $y'\in B_{\theta}(\bar{x}')$,
\begin{equation}\label{S2-HP}
    |[u(y',\, \cdot) \leq M] \cap B_{\rho}(\bar{x}'')| \leq (1-\alpha) |B_{\rho}|,\quad \text{being} \quad \theta= \rho^{\frac{p}{2}} (\delta M)^{\frac{2-p}{2}},
\end{equation}\noindent and $Q_{2\theta,4\rho}(\bar{x}) \subset \Omega$. Then for every $\nu \in (0,1)$ there exist numbers $K>1$ and $\delta_o\in (0,1)$ depending only on the data and $\alpha,\nu,n$ such that either
\begin{equation}\label{alternativeSH2}
    M\leq K \rho, \quad \text{or}
\end{equation}
\begin{equation}\label{S2-TH}
    |[u(x',\cdot) \leq \delta_o M] \cap B_{2\rho}(\bar{x}'')| \leq \nu |B_{2\rho}(\bar{x}'')|, \quad \quad \forall_{ae}\,\,  x' \in B_{\theta}(\bar{x}').
\end{equation}\noindent 

\end{lemma}

\begin{proof} We divide the proof into three steps. \vskip0.2cm \noindent 
{\small{\it STEP 1. NORMALIZATION AND LOGARITHMIC ESTIMATE.}}
\vskip0.2cm \noindent 
Let us introduce the change of variables $\Phi: Q_{2\theta,4\rho}(\bar{x}) \rightarrow Q_{1,1}(0)$ given by
\[
x' \rightarrow \frac{x'-\bar{x}'}{2\theta}, \quad \quad x'' \rightarrow \frac{x''-\bar{x}''}{4\rho}, \quad \quad v= \frac{u}{M}.\] 
The new function $v$ satisfies the following equation in $Q_{1,1}$,
\begin{equation}\label{NE1}
\sum_{i=1}^s \partial_{ii} v + \sum_{i=s+1}^N \partial_i ( \tilde{A}_i(x,v,\nabla v))=0,
\end{equation}
with structure conditions 
\begin{equation} \label{E2-structure}
\begin{cases}
 \sum_{i=s+1}^N \tilde{A}_i(x,v,\nabla v) \, \partial_i v \ge \tilde{C}_1 \sum_{i=s+1}^N |\partial_i v|^p-\tilde{C} ( \rho/M )^{p},\\
 |A_i(x,v,\nabla v)| \leq \tilde{C}_2 |\partial_i v|^{p-1}+\tilde{C} (  \rho/M )^{p-1}.
\end{cases} 
\end{equation} \noindent
By transformation $\Phi$ inequality \eqref{S2-HP} turns into
% \begin{equation} \label{ciccio0}
%  |[v(x',\cdot) \leq 1] \cap B_1(0'')| \leq (1-\alpha) |B_1(0'')|, \quad \quad \forall_{\text{ae}}\,\,   x' \in B_{1}(0'),
% \end{equation} \noindent that is easily seen to imply the measure estimate
\begin{equation}\label{ciccio0}
   |[v(x',\,\cdot) > 1] \cap B_{1/4}(0'')| > \alpha |B_{1/4}(0'')|, \quad \quad \forall_{\text{ae}} \, \, x' \in B_{1/2}(0').
\end{equation} \noindent The expansion of positivity relies on the following simple fact. The inequality above implies the measure estimate for $x'$ in $B_{1/2}(0')$
\begin{equation}\label{ciccio}
   |[v(x',\,\cdot) > 1] \cap B_{1}(0'')| > (w_{N,s} \,4^{(N-s)})^{-1}\alpha |B_{1}(0'')|=:\tilde{\alpha} |B_{1}(0'')|.
\end{equation}\noindent 
Let now $\zeta(x'')\in C_o^{\infty}(B_1(0''))$ be a convex cut-off function between balls $B_1$ and $B_{1/2}$, i.e.
\[0 \leq \zeta \leq 1, \quad \quad \zeta_{|B_{1/2}(0'')} \equiv 1, \quad \quad  |\partial_i \zeta| \leq 2 \quad \forall i \in \{s+1,...,N\}.\]
Let us fix numbers $j^* \in \N$ and $\epsilon \in (0,1)$ and for $j=1,2,..,j^*$ we define 
\[
f((v-\e^j)_{-})=[(\e^j(1+\e)-(v-\e^j)_{-}]^{1-p}.
% \quad \quad \text{therefore assigning levels} \quad k=\e^j.
\]
% \[
% 0<\frac{d}{d(v-\e^j )_{-}} \,f((v-\e^j)_{-})=(p-1)[(\e^j(1+\e)-(v-\e^j)_{-}]^{-p} \, \chi_{[v<\e^j]}(x) \leq \e^{-j^+ p}. \]\vskip0.4cm \noindent
Then, inequality \eqref{energy2} reads, for $t \in (0,1/2)$,
\begin{equation}\label{ciccio2}
\begin{aligned}
\gamma^{-1}& \sum_{i=s+1}^N  \int \int_{Q_{t,1}} \frac{|\partial_i (v-\e^j)_{-}|^{p}}{[\e^j(1+\e)-(v-\e^j)_{-}]^p} (t^2-|x'|^2) \zeta^p(x'') dx \leq\\
& t^{s+1}\frac{d}{dt} \bigg[t^{-s}\int \int_{Q_{t,1}}  \bigg(\int_0^{(v-\e^j)_{-}} f(\tau) d\tau\bigg) \zeta(x'')^p\, dx  \bigg]+\\
& \int \int_{Q_{t,1}} \bigg{\{}\sum_{i=s+1}^N (t^2-|x'|^2) ( |\partial_i \zeta(x'')|^p+\zeta(x'')^p )+\tilde{C} \bigg(\frac{\rho}{M\e^{j^*}}\bigg)^{p} (t^2-|x'|^2) \bigg{\}} dx.
\end{aligned}
\end{equation} \noindent Last two terms on the right of \eqref{ciccio2} can be reduced, by assuming $M\e^{j*} >\rho$, to
\[
\sum_{i=s+1}^N t^{-s-1} |Q_{t,1}| t^2+ \tilde{C} \bigg(\frac{\rho}{M \e^{j*}} \bigg)^p t^{-s-1} |Q_{t,1}| t^2 \leq \gamma t.
\] We observe that first integrands on the left of \eqref{ciccio2} are the directional derivatives of
\[
g= \bigg[ \ln \bigg( \frac{(1+\e)\e^j}{(1+\e)\e^j-(v-\e^j)_{-}} \bigg) \bigg]^p,
\] and, as $g \in W^{1,p}(B_{1}(0''))$ vanishes in $[v>e^{j}] \cap B_1|$, we can apply the weighted Poincar\'e inequality \eqref{WP}, using \eqref{ciccio} to estimate the term $|\mathcal{E}|=|[v>\e^j]\cap B_{1}|\ge|[v>1]\cap B_{1}|$.\newline \noindent A precise analysis of this last simple fact reveals its correspondence with \eqref{ciccio0}-\eqref{ciccio} in terms of expansion of positivity. Putting all the pieces of the puzzle into \eqref{ciccio2} we arrive finally to the following logaritmic estimate, valid for each $j=1,..,j^*$, and $0<t<\theta$,
\begin{equation}\label{ciccio3}
\begin{aligned}
    t^{-s-1} \int \int_{Q_{t,1}}& \ln^p \bigg( \frac{(1+\e)\e^j}{(1+\e)\e^j-(v-\e^j)_{-}} \bigg) (t^2-|x'|^2) \zeta^p(x'') dx \leq\\
    &\gamma\frac{d}{dt} \bigg( \int \int_{Q_{t,1}}t^{-s} \int_0^{(v-\e^j)_{-}} \frac{d \tau}{[(1+\e)\e^j-\tau]^{p-1}}  \zeta^p(x'')dx \bigg) + \gamma t.
    \end{aligned}
\end{equation}  

\vskip0.2cm \noindent 
{\small{\it STEP 2. FIRST ALTERNATIVE.}}
\vskip0.2cm

\noindent 
Let us define $A_j(t)= [v<\e^j] \cap Q_{t,1}(0)$ and let 
\begin{equation}\label{Y}
    y_j= \sup_{0<t<1/2} Y_j(t), \quad \text{where} \quad Y_j(t)= t^{-s} \int \int_{A_j(t)} \zeta^p(x'') \, dx.
\end{equation} \noindent 
Now we show that if $M \ge K \rho$, then there exists a number $\xi=\xi(\nu) \in (0,1)$ such that $y_{j+1} \leq \max \{\nu, (1-\xi) y_j  \}$ for each $j=1,..,j^*$. This, with a standard iteration procedure will end the proof. So we proceed by assuming $ y_{j+1} > \nu$, and by continuity of the integral we can choose $t_o \in (0,1/2)$ such that $Y_{j+1}(t_o)=y_{j+1}$ and divide the argument in two alternatives. Let
\begin{equation}\label{alt1}
        \Psi(t):= t^{-s} \int \int_{Q_{t,1}} \bigg( \int_0^{(v-\e^{j})_{-}} \frac{d\tau}{[(1+\e)\e^{j}-\tau]^{p-1}} \bigg) \zeta^p(x'') dx,
    \end{equation} \noindent and suppose $\Psi'(t_o) \leq 0$. Then inequality \eqref{ciccio3} implies that for each fixed $\sigma \in (0,1)$ 
\begin{equation} \begin{aligned} \label{ciccio4}
&t_o^{-s-1}t_o^2[1-(1-\sigma)^2]\ln^p \bigg( \frac{1+\e}{2\e} \bigg) \int \int_{A_{j+1}(1-\sigma)t_o}  \zeta^p(x'') \, dx \leq\\
&t_o^{-s} \int \int_{A_{j}(1-\sigma)t_o} \ln^p \bigg( \frac{(1+\e)\e^{j}}{(1+\e)\e^{j}-(v-\e^{j})_{-}} \bigg) \zeta^p(x'') (t_o^2-|x'|^2))\, dx \leq \gamma t_o,
\end{aligned} \end{equation} \noindent giving the estimate
\begin{equation}
    \begin{aligned}
    t_o^{-s} \int \int_{A_{j+1}(1-\sigma)t_o} \zeta^p(x'') \, dx \leq \gamma \sigma^{-2} \ln^{-p} \bigg( \frac{1+\e}{2\e} \bigg).
    \end{aligned}
\end{equation} \noindent 
Now we determine $\sigma, \e$ small enough to get $y_{j+1} \leq \nu$. Indeed, by $|Q_{t,1}|\leq \gamma t^{-s} |B_1(0'')|$ we see that the following estimate holds
\begin{equation}\label{ciccio5}
    \begin{aligned}
    y_{j+1}&=Y_{j+1}(t_o)= t_o^{-s} \int \int_{A_{j+1}(t_o)} \zeta^p(x'') \, dx \leq\\&
    t_o^{-s} \int \int_{Q_{(1-\sigma)t_o,1}} \chi_{[v<\e^{j+1}]}\zeta^p(x'') \, dx+ t_o^{-s} |Q_{t_o,1}/Q_{(1-\sigma)t_o,1}| \leq\\
    &\gamma \bigg{\{} \sigma^{-2}\ln^{-p} \bigg( \frac{1+\e}{2\e} \bigg)+ \sigma (N-s) |B_1(0'')|\bigg{\}} \leq \nu,
    \end{aligned}
\end{equation} \noindent for the choices
\begin{equation}\label{sigma-epsilon}
\sigma= \frac{\nu}{2(N-s) \gamma}, \quad \quad \e= ( e^{\frac{2\gamma}{\sigma^2 \nu p}}-1)^{-1}.
\end{equation} \noindent

\vskip0.2cm \noindent 
{\small{\it STEP 3. SECOND ALTERNATIVE.}}
\vskip0.2cm

\noindent 
Let us suppose now that $\Psi'(t_o) >0$ and 
% \begin{equation}\label{alt2}
% \Psi(t_o)= \frac{d}{dt}_{|t=t_o} \bigg[t^{-s} \int \int_{Q_{t,1}} \bigg( \int_0^{(v-\e^j)_{-}} \frac{d\tau}{[(1+\e)\e^j-(v-\e^j)_{-}]^{p-1}} \bigg) \zeta^p(x'') dx  \bigg] > 0.
% \end{equation}\noindent 
that there exists $t_*= \inf \{ t\in(t_o,1/2) | \quad \Psi'(t) \leq 0   \}$, so that by definition $\Psi$ is monotone increasing before $t_*$ and we have
\begin{equation}\label{ciccio6}
\begin{aligned}
   t_o^{-s} \int \int_{Q_{t_o,1}}& \bigg( \int_0^{(v-\e^{j})_{-}} \frac{d\tau}{[(1+\e)\e^{j}-\tau]^{p-1}} \bigg) \zeta^p(x'') dx \leq\\
   &t_*^{-s} \int \int_{Q_{t_*,1}} \bigg( \int_0^{(v-\e^{j})_{-}} \frac{d\tau}{[(1+\e)\e^{j}-\tau]^{p-1}} \bigg) \zeta^p(x'') dx .
   \end{aligned}
\end{equation}\noindent For the time $t_*$ similar estimates to \eqref{ciccio3}, \eqref{ciccio4}, \eqref{ciccio5} hold and we obtain similarly that 
\begin{equation}\label{ciccio7}
    t_*^{-s}\int \int_{A_{j}(t_*)} \zeta^p(x'') \chi_{[(\e^{j}-v)>\e^{j}\tau]} dx \leq \nu/4+ \gamma \nu^{-2} \ln^{-p} \bigg(\frac{1+\e}{1+\e-\tau}  \bigg) \leq \nu/2.
\end{equation} \noindent In this case the value of $\e$ has been already chosen, so that last inequality is valid provided we restrict to levels $\tau$ such that
\begin{equation}
    \gamma \nu^{-2} \ln^{-p} \bigg(\frac{1+\e}{1+\e-\tau}  \bigg) \leq \nu/4 \quad \iff \quad \tau \ge (1+\e) [1-e^{-h(\nu)}]=:\bar{\tau},
\end{equation}\noindent for a function $h(\nu)= o(\nu^{3/p})$. Finally, we use \eqref{ciccio7} and \eqref{ciccio6} together with Fubini's theorem and a change of variables to estimate
\begin{equation} \label{ciccione}
    \begin{aligned}
    &t_*^{-s} \int \int_{A_{j}(t_*)} \bigg( \int_0^{(v-\e^{j})_{-}} \frac{d\tau}{[(1+\e)\e^{j}-\tau]^{p-1}} \bigg) \zeta^p(x'') dx=\\
    % &    t_*^{-s}  \int \int_{A_{j+1}(t_*)} \zeta^p(x'') \bigg{\{} \int_0^{\e^{j+1}} \frac{\chi_{[\e^{j+1}-v \ge  \tau]}\, d \tau}{[(1+\e)\e^{j+1}-\tau]^{p-1}} \bigg{\}} dx=\\
    % &\int_0^{\e^{j+1}} \frac{1}{[(1+\e)\e^{j+1}-\tau]^{p-1}} \bigg{\{} t_*^{-s} \int \int_{A_{j+1}(t_*)} \zeta^p(x'') \chi_{[\e^{j+1}-v\ge  \tau]} dx \bigg{\}} d\tau=\\
    & \e^{j(2-p)} \int_0^1 \frac{1}{[1+\e-\tau]^{p-1}}  \bigg{\{} t_*^{-s} \int \int_{A_{j}(t_*)} \zeta^p(x'') \chi_{[\e^{j}-v\ge  \e^{j}\tau]} dx \bigg{\}} d\tau \leq\\
    & \e^{j(2-p)}   \bigg[y_j \int_0^{\bar{\tau}} \frac{d\tau}{[1+\e-\tau]^{p-1}}  + \frac{\nu}{2} \int_{\bar{\tau}}^{1} \frac{d\tau}{[1+\e-\tau]^{p-1}}\bigg] \leq\\
    &\e^{j(2-p)}   y_j \bigg[\int_0^{\bar{\tau}} \frac{d\tau}{[1+\e-\tau]^{p-1}}  + \frac{1}{2}\int_{\bar{\tau}}^{1} \frac{d\tau}{[1+\e-\tau]^{p-1}}\bigg],
    \end{aligned}
\end{equation} \noindent where we used $y_{j}> Y_{j+1}(t_*)$ in first inequality and $y_j>y_{j+1}>\nu$ in last one. By an easy manipulation we see that we can bound from below \eqref{ciccio6}, indeed
\begin{equation}\label{manipulation}
    \begin{aligned}
    &\int_0^{(v-\e^{j})_{-}} \frac{d\tau}{[(1+\e)\e^{j}-\tau]^{p-1}}\ge
    % & \e^{(j+1)(2-p)} \int_0^{1-\e} \frac{d\tau}{[1+\e-\tau]^{p-1}}= \e^{(j+1)(2-p)} \bigg(\int_0^{1} \frac{d\tau}{[\e+\tau]^{p-1}} -\int_0^{\e}\frac{d\tau}{[\e+\tau]^{p-1}} \bigg) = \\
    \e^{j(2-p)} (1-\e^{2-p}) \int_0^1\frac{d\tau}{[\e+\tau]^{p-1}}\quad .
    \end{aligned}
\end{equation} \noindent So we gather \eqref{ciccio6}, \eqref{ciccione} and \eqref{manipulation} together to get
\begin{equation}\label{final}
    \begin{aligned}
    y_{j+1}& \e^{j(2-p)} (1-\e^{2-p}) \int_0^1\frac{d\tau}{[\e+\tau]^{p-1}} \leq\\
    % &t_o^{-s} \int \int_{Q_{t_o,1}} \zeta^p(x'') \bigg( \int_0^{(v-\e^{j+1})_{-}} \frac{d\tau}{[(1+\e)\e^{j+1}-\tau]^{p-1}} \bigg) dx\leq\\
    &\e^{j(2-p)}y_j \bigg(\int_0^1 \frac{d\tau}{(\e+\tau)^{p-1}}-\frac{1}{2}\int_0^{1-\bar{\tau}} \frac{d\tau}{(\e+\tau)^{p-1}}  \bigg) \leq \\
    &\e^{j(2-p)}y_j \bigg(\int_0^1 \frac{d\tau}{(\e+\tau)^{p-1}}-\frac{(1-\bar{\tau})}{2}\int_0^{1} \frac{d\tau}{(\e+\tau(1-\bar{\tau}))^{p-1}}  \bigg) \leq\\
    &\e^{j(2-p)}y_j (1-(1-\bar{\tau})/2)  \int_0^1 \frac{d\tau}{(\e+\tau)^{p-1}}.
    \end{aligned}
\end{equation} \noindent 
% having used for the term on the right the simple estimate
% \[
% \int_0^{1-\bar{\tau}} \frac{d\tau}{(\e+\tau(1-\bar{\tau}))^{p-1}} = \int_0^1 \frac{(1-\bar{\tau})d\tau}{(\e+(1-\bar{\tau})\tau)^{p-1}} \ge (1-\bar{\tau}) \int_0^1 \frac{d\tau}{(\e+\tau)^{p-1}}.
% \]  
\noindent 
Finally, this implies
\[
y_{j+1} \leq \bigg(\frac{1-(1-\bar{\tau})/2}{1-\e^{2-p}}  \bigg) y_j=:(1-\xi) y_j, \quad\quad  \text{with} \quad \quad \xi =1-\bigg[\frac{1-(1-\bar{\tau})/2}{1-\e^{2-p}}  \bigg]<1,
\] redefining $\e = \min \{1/2, \e\}$ if needed. We prove now that if $t_*$ does not exist, then the iteration inequality above is still satisfied. Indeed in case no such $t_{*}$ exists then $\Psi'(t)>0$ for all $t \in [t_o,1/2]$ and therefore $\Psi(t_o) \leq \Psi(1/2)$. Moreover by simple calculations analogous to \eqref{manipulation}, \eqref{final} we recover the estimates
\begin{equation}\label{mick}
    \Psi(t_o) \ge \epsilon^{j(2-p)} (1-\epsilon^{2-p}) \bigg(\int_0^1 \frac{d\tau}{(\epsilon+\tau)^{p-1}}\bigg) y_{j+1},
\end{equation}
\begin{equation}\label{jagger}
    \Psi(1/2) \leq 2^s \epsilon^{j(2-p)} \int_0^1 \frac{d\tau}{(\epsilon+\tau)^{p-1}} \int \int_{Q_{1/2,1}} \chi_{[v \leq \epsilon^{j}]}(x) \zeta^p(x'') \, dx.
\end{equation} \noindent So, by combining $\Psi(t_o) \leq \Psi(1/2)$ with \eqref{mick}, \eqref{jagger} we obtain the inequality
\begin{equation}\label{will}
    y_{j+1} \leq \frac{2^s}{1-\epsilon^{2-p}} \int \int_{Q_{1/2,1}} \chi_{[v\leq \epsilon^{j}]}\zeta^p(x'')\, dx \leq \frac{2^s}{1-\epsilon^{2-p}} \int \int_{Q_{1/2,1}} \chi_{[v\leq \epsilon]}\zeta^p(x'')\, dx.
\end{equation} \noindent If we test equation \eqref{mon-energy} with $\psi= [\epsilon(1+\epsilon)-(v-\epsilon)_{-}]^{1-p}\zeta^p(x'')(1-|x'|^2)_+^2$ and use Young inequality we can derive, similarly to \eqref{ciccio3}, the estimate
\begin{equation*}\label{smith}
\begin{aligned}
&\int \int_{Q_{1/2,1}} \ln^p \bigg(\frac{1+\epsilon}{\epsilon(1+\epsilon)-(v-\epsilon)_{-}}  \bigg)\zeta^p(x'')\,dx\leq \\
&\gamma \sum_{i=s+1}^N \int \int_{Q_{1,1}} \bigg(\frac{|\partial_i (v-\epsilon)_{-}| }{[\epsilon(1+\epsilon)-(v-\epsilon)_{-}]}\bigg)^p\zeta^p(x'') (1-|x'|^2)\, dx \leq\\
&\gamma \int \int_{Q_{1,1}} \bigg{\{}[\epsilon(1+\epsilon)-(v-\epsilon)_{-}]^{2-p} \zeta^p(x'')+(1-|x'|^2) \bigg{\}} dx+ \tilde{C} \bigg( \frac{\rho}{M \epsilon^{j^{*}}} \bigg)^p |Q_{1,1}| \leq \gamma.
\end{aligned}
\end{equation*}\noindent From this, by \eqref{jagger} we obtain
\[
\int\int_{Q_{1/2,1}} \chi_{[v\leq \epsilon]} \zeta^p(x'') \, dx \leq \gamma \ln^{-p} \bigg(\frac{1+\epsilon}{2\epsilon} \bigg),
\] and therefore estimating \eqref{will} from above we conclude that $y_{j+1} \leq \nu$ by choosing $\epsilon$ small enough. Finally, for the sake of readability we just remark that in case $\Psi$ is not regular enough it is possible to perform the same argument above by substituting $\Psi'$ with its right Dini derivative, as in (\cite{DB-Chen}, Sec. 7). 

\vskip0.2cm \noindent 
{\small{\it CONCLUSION.}}
\vskip0.2cm

\noindent Both the alternatives imply the estimate
\[
y_{j+1} \leq \max\{ \nu,(1-\xi) y_j \}, \quad \forall j=1,..,j^*, \quad \xi=\xi(\nu) \in (0,1).
\] Iterating this inequality we arrive at
\[
y_{j^*} \leq \max\{ \nu, (1-\xi)^{j^*-1} y_1 \} 
\] and since $y_1\leq 1$, we choose $j^*$ such that $(1-\xi)^{j^*-1}\leq \nu$ to get for each $x'\in (0,1/2)$ the estimate
\begin{equation} \label{A}
|[v(x',\cdot)<\e^{j^*}] \cap B_1(0'')| \leq y_{j^*} \leq \nu |B_1(0'')|.
\end{equation} \noindent The inverse transformation $\Phi^{-1}$ turns the obtained estimate into \eqref{S2-TH} with $\delta_o=\e^{j^*}$, therefore finishing the proof of Lemma \ref{SH2-Lemma}.

\end{proof}

\subsection{Conclusion of the Proof of Theorem \ref{shrinkingTHM}}
Let $\bar{x} \in \Omega$, $\rho>0$ and $\alpha \in (0,1)$. We suppose that for $\delta(\alpha) \in (0,1)$ and $M >0$ we have the information
\[
Q_{2\theta,4\rho}(\bar{x}) \subset \Omega, \quad \text{\&} \quad |[u \leq M] \cap Q_{\theta,\rho}(\bar{x})| \leq (1-\beta(\alpha)) |Q_{\theta,\rho}|, \quad \text{being} \quad \theta= \rho^{\frac{p}{2}} (\delta M)^{\frac{2-p}{p}}.
\] By Lemma \ref{SH1-Lemma} there exist numbers $s_1, \delta_1>0$ depending only from the data and $\alpha$ such that either $M \leq \rho$ or 
\[
|[u(y',\cdot) \leq 2^{-s_1} M] \cap B_{\rho}(\bar{x}'')| \leq (1-\alpha^4/2) |B_{\rho}(\bar{x}'')|, \quad \forall_{\text{ae}}\,\, y' \in B_{(1-\delta_1)\theta} (\bar{x}').
\] For $\nu$ as in \eqref{nu}, we apply Lemma \ref{SH2-Lemma} with
\[
\bar{M}= 2^{-s_1}M, \quad \beta= \beta{\alpha},\quad \bar{\delta} = (1-\delta_1)^{\frac{2}{p-2}},\quad  \bar{\theta}= \rho^{\frac{p}{2}} (\bar{\delta}M)^{\frac{2-p}{2}}= (1-\delta_1)\theta,
\] so that there exist numbers $K>1$, $ \delta_o \in (0,1)$ depending on the data and $\alpha,\nu$ such that either $M\leq K \rho$ or 
\[
|[u(x', \cdot) \leq \delta_o M] \cap B_{4\rho}(\bar{x}'')| \leq \nu |B_{4\rho} (\bar{x}'')|, \quad \forall_{ae}\, \, x' \in B_{\bar{\theta}}(\bar{x}').
\] \noindent To recover the correct intrinsic geometry (see Remark \ref{R-geometry}), we cut the slice-wise information on polydisc $Q_{\bar{\theta},4\rho}(\bar{x})$ to an information on a polydisc which is smaller along the nondegenerate variables, by
\[
\theta_o^2= (4\rho)^{p} (\bar{\delta}_o M)^{2-p}, \quad \bar{\delta}_o=4^{\frac{-p}{2-p}}\delta_o, 
\] so that $\theta_o \leq \bar{\theta}$, increasing $j^*$ in \eqref{A} in case of need. Finally we apply the {\it Critical Mass} Lemma \ref{DG} to end the proof.

\begin{center}
\section{H\"older Continuity. Proof of Theorem \ref{holderTHM}\label{4}}
\end{center} \noindent 
We begin with the accommodation of degeneracy. Let $x_0 \in \Omega$ be an arbitrary point, $M= \sup_{\Omega} |u|$ and $\rho>0$ such that \[
Q_{\rho}[M] (x_0):= Q_{\rho^{\frac{p}{2}}(2M)^{\frac{2-p}{2}},\rho}(x_0) \subset \Omega.
\] Set
\[
\mu_+:= \sup_{\Omega} u, \quad \mu_{-}:= \inf_{\Omega} u\, \quad \omega= \mu_+-\mu_-.
\] Now for $\alpha=1/2$ we fix a number $\delta$ as defined in Theorem \ref{shrinkingTHM}, let $\theta= \rho^{\frac{p}{2}} (\delta \omega)^{\frac{2-p}{2}}$, and consider the following two alternatives:
\begin{equation}\label{HC-alt1}
    |[u\leq \mu_-+ \omega/2] \cap Q_{\theta,\rho}(x_0)| \leq \frac{1}{2}|Q_{\theta, \rho}|, 
\end{equation}\noindent or
\begin{equation}\label{HC-alt2}
    |[u\ge \mu_+ - \omega/2] \cap Q_{\theta,\rho}(x_0)| \leq \frac{1}{2}|Q_{\theta, \rho}|.
\end{equation} \noindent 
From this, in both cases we may apply Theorem \ref{shrinkingTHM}: in case of measure estimate \eqref{HC-alt1} we apply it to the function $v^+=(u-\mu_-)$, while in case of measure estimate \eqref{HC-alt2} $v^-=(\mu_+-u)$. Both these functions are solutions to an equation similar to \eqref{E1}-\eqref{E1-structure}, and therefore the aforementioned Theorem can be applied, implying a reduction of oscillation 
\[
\osc_{Q_{\eta,\rho/4}} u \leq \bigg(1-\frac{\delta_o}{4} \bigg) \omega:= \delta_* \omega, \quad \quad \eta^2=(\rho/4)^p (\delta \omega)^{2-p}.
\] Once we have this kind of controlled reduction of oscillation the whole procedure can be iterated in nested shrinking polydiscs and the rest of the proof is standard (see for instance Theorem 3.1 in \cite{DB2}, Chap. X).
\vskip0.4cm \noindent 
{\centering{\section{$L^1-L^{\infty}$ estimates. Proof of Theorem \ref{l1-linftyTHM}\label{5}}}}
\noindent Let us fix a point $\bar{x} \in \Omega$ and numbers $\rho, \theta>0$ such that $ Q_{8\theta,8\rho}(\bar{x}) \subset \Omega$. Let $y' \in B_{\theta/2} (\bar{x}')$ be a Lebesgue point for the function \[ y \rightarrow \int_{B_r(\bar{x}'')} u(y,x'') \zeta(x'') \, dx''\, .\] For $\sigma \in (0,1)$ we consider a generic radius $r$ such that $\rho \leq (1-\sigma)r \leq 2 \rho$. Let $\zeta(x'') \in C_o^{\infty}(B_r(\bar{x}))$, $0\leq \zeta \leq 1$, $\zeta \equiv 1$ in $B_{(1-\sigma)r}(\bar{x}'')$ be a cut off function relative to the last $N-s$ variables between $B_{(1-\sigma)r}$ and $B_r$ satisfying
\[
|\partial_i \zeta |\leq \frac{\gamma}{\sigma r} , \quad \forall i\in \{s+1,..,N\}.
\] We divide the proof in three steps. For ease of notation, let us call $\Q_{t,r}:= Q_{t,r}(y',\bar{x}'')$ and 
\[\eta= \bigg(\frac{\theta^2}{\rho^p}\bigg)^{\frac{1}{2-p}}.\]
 \noindent 
{\small{\it STEP 1. AN INTEGRAL INEQUALITY}}
\vskip0.2cm \noindent 
We test the equation \eqref{E1}-\eqref{E1-structure} with $\phi(x)= (t^2-|x'-y'|^2)_+\zeta^p(x'')$ for $0<t<2\theta$ and use Green's formula \eqref{Stokes} with $f(u) \equiv 1$ to get
\begin{equation}\label{passo1}
\begin{aligned}
t^{s+1} \frac{d}{dt} \bigg(|B_t(0')|^{-1} \int  \int_{\Q_{t,r}} u \zeta^p \, dx \bigg) &=  \sum_{i=1}^s \int \int_{\Q_{t,r}} (\partial_i u)  (x^{'}_i-y^{'}_i)  \zeta^p  \, dx\leq \\
    & \frac{p}{2}\sum_{s+1}^N \int \int_{\Q_{t,r}} (t^2-|x'-y'|^2)_+\, A_i(x, \nabla u)\zeta^{p-1} \partial_i \zeta \, dx.
    \end{aligned}
\end{equation} \noindent Now we estimate \eqref{passo1} by use of structure conditions \eqref{E1-structure} to get
\begin{equation}\label{passo3}
    \begin{aligned}
    &t^{s+1} \frac{d}{dt} \bigg(t^{-s} \int \int_{\Q_{t,r}} u \zeta^p \, dx \bigg)\leq \bigg( \int \int_{\Q_{t,r}} (u+\eta)^{\beta(p-1)} \, dx \bigg)^{\frac{1}{p}} \times\\
    &\times \frac{\gamma(\omega_s)}{\sigma \rho} \sum_{i=s+1}^N \bigg( \int \int_{\Q_{t,r}}(u+\eta)^{-\beta} |\partial_i u|^p (t^2-|x'-y'|^2)_+^{\frac{p}{p-1}} \zeta^p\, dx\bigg)^{\frac{p-1}{p}}+\frac{\gamma t^2}{\sigma \rho} |\Q_{t,r}|,
    \end{aligned}
\end{equation} \noindent using H\"older inequality and multiplying and dividing for $(u+\eta)^{\beta (p-1)/p}$, $\beta>0$.\vskip0.2cm \noindent Let us estimate the second integral term in \eqref{passo3}. Let $\zeta(x'')$ be as before and let us test the equation \eqref{E1}-\eqref{E1-structure} with 
\[
\psi(x)= (u(x)+\eta)^{1-\beta} (t^2-|x'-y'|^2)_+^{\frac{p}{p-1}} \zeta^p(x'') \in W^{1,\bf{p}}_o(\Q_{t,r}).
\] Using structure conditions \eqref{E1-structure} and Young's inequality, we get
\begin{equation}\label{passo4}
    \begin{aligned}
     I:=&\gamma^{-1}\int \int_{\Q_{t,r}} (u+\eta)^{-\beta} (t^2-|x'-y'|^2)_+^{\frac{p}{p-1}} \bigg{\{} \sum_{i=1}^s|\partial_i u|^2 + \sum_{i=s+1}^N |\partial_i u|^p \bigg{\}} \, dx \leq \\
     &  \int \int_{\Q_{t,r}} (u+\eta)^{-\beta} (t^2-|x'-y'|^2)_+^{\frac{p}{p-1}} \zeta^p\, dx +\\
     &+\sum_{i=1}^s \int \int_{\Q_{t,r}}(\partial_i u)(u+\eta)^{1-\beta}(t^2-|x'-y'|^2)_+^{\frac{p}{p-1}-1} (x^{'}_i-y^{'}_i) \zeta^p \, dx +\\
     &+\sum_{i=s+1}^N\frac{1}{\sigma \rho} \int \int_{\Q_{t,r}} (|\partial_i u|^{p-1}+1)(u+\eta)^{1-\beta} (t^2-|x'-y'|^2)_+^{\frac{p}{p-1}} \zeta^{p-1}\, dx.\end{aligned}
\end{equation} \noindent 
Applying repeatedly Young's inequality and reabsorbing on the left the terms involving energy estimates we obtain
\begin{equation}\label{passo5}
    \begin{aligned}
    I \leq &  \int \int_{\Q_{t,r}} (u+\eta)^{2-\beta} (t^2-|x'-y'|^2)_+^{\frac{p}{p-1}-2} |x'-y'|^2 \zeta^p \, dx + \\
    &+\frac{1}{(\sigma \rho)^p} \int \int_{\Q_{t,r}} (u+\eta)^{p-\beta} (t^2-|x'-y'|^2)_+^{\frac{p}{p-1}}  \zeta^p \, dx +\\
    &+  \int \int_{\Q_{t,r}} (u+\eta)^{-\beta} (t^2-|x'-y'|^2)_+^{\frac{p}{p-1}} \zeta^p\, dx =\\
    &= I_1+I_2+\gamma \eta^{-\beta} t^{\frac{2p}{p-1}}|\Q_{t,r}|.
    \end{aligned}
\end{equation} \noindent We estimate separately the various terms. For the first one we have
\begin{equation*}
    \begin{aligned}
    I_1=  \int \int_{\Q_{t,r}} (u+\eta)^{2-\beta} (t^2-|x'-y'|^2)_+^{\frac{p}{p-1}-2} |x'-y'|^2 \zeta^p \, dx \leq \gamma t^{\frac{2}{p-1}} \int \int_{\Q_{t,r}} (u+\eta)^{2-\beta} \, dx.
    \end{aligned}
\end{equation*} \noindent Next we use $1<p<2$ to split $(u+\eta)^{p-\beta}=(u+\eta)^{p-2} (u+\eta)^{2-\beta}\leq \eta^{p-2}(u+\eta)^{2-\beta}$ to get
\begin{equation*} 
    \begin{aligned}
    I_2 \leq\frac{\gamma }{(\sigma \rho)^p}\eta^{p-2} t^{2\frac{p}{p-1}} \int \int_{\Q_{t,r}} (u+\eta)^{2-\beta}\, dx \leq \frac{\gamma}{\sigma^p} t^{\frac{2}{p-1}}\int_{\Q_{t,r}} (u+\eta)^{2-\beta}\, dx. 
    \end{aligned}
\end{equation*} \noindent Inequalities above about $I_1,I_2$ and \eqref{passo5} lead us to the formula
\begin{equation}\label{passo6}
\begin{aligned}
\sum_{i=s+1}^N \int \int_{\Q_{t,r}} &|\partial_i u|^p (u+\eta)^{-\beta} (t^2-|x'-y'|^2)_+^{\frac{p}{p-1}} \zeta^p \, dx \leq\\
& \gamma \sigma^{-p} t^{\frac{2}{p-1}} \int \int_{\Q_{t,r}} (u+\eta)^{2-\beta} \, dx + \gamma \sigma^{-p} \eta^{-\beta} t^{\frac{2p}{p-1}}|\Q_{t,r}|.
\end{aligned}
\end{equation}\noindent 
We put \eqref{passo6} inside \eqref{passo3}, summing $(N-s)$ times the same quantity, to obtain 
\begin{equation}\label{passo7}
    \begin{aligned}
    &t^{s+1} \frac{d}{dt} \bigg(t^{-s} \int \int_{\Q_{t,r}} u \zeta^p \, dx \bigg)\leq \, \frac{\gamma t^2}{\sigma \rho} |\Q_{t,r}|+\\
    &+\frac{\gamma}{\sigma^{p} \rho} \bigg( t^{\frac{2}{p-1}} \int \int_{\Q_{t,r}} (u+\eta)^{2-\beta} \, dx + \eta^{-\beta} t^{\frac{2p}{p-1}}|\Q_{t,r}| \bigg)^{\frac{p-1}{p}} \bigg( \int \int_{\Q_{t,r}} (u+\eta)^{\beta(p-1)} \, dx \bigg)^{\frac{1}{p}}.
    \end{aligned}
\end{equation} \noindent 
Now we evaluate by Jensen's inequality separately for $\alpha= 2-\beta,\,\beta(p-1)$, the terms
\begin{equation*}\label{passo8}
    \begin{aligned}
    \int \int_{\Q_{t,r}}& (u+\eta)^{\alpha} \, dx   \leq|\Q_{t,r}| \bigg( \dashint \dashint_{\Q_{t,r}} (u+\eta)\, dx \bigg)^{\alpha}.
    \end{aligned} \end{equation*}
    
%     \begin{equation*}\label{passo9}
%     \begin{aligned}
%     \int & \int_{\Q_{t,r}} (u+\eta)^{\beta(p-1)} \, dx 
%   \leq|\Q_{t,r}| \bigg( \dashint \dashint_{\Q_{t,r}} (u+\eta)\, dx \bigg)^{\beta(p-1)}.
%     \end{aligned} \end{equation*}\noindent 
\noindent Let 
\[
\A:=\bigg( \dashint \dashint_{\Q_{t,r}} (u+\eta) \, dx\bigg)= \bigg(  \dashint \dashint_{\Q_{t,r}} u\, dx \bigg)+ \eta,
\] and therefore we estimate from above inequality \eqref{passo7} by
\begin{equation*} \label{passo10}
    \begin{aligned}
    t^{s+1} \frac{d}{dt}  \bigg(t^{-s} \int \int_{\Q_{t,r}} u \zeta^p \, dx \bigg)\leq  \frac{\gamma}{\rho \sigma^p} |\Q_{t,r}| \bigg{\{} t^{\frac{2}{p}} \A^{2(\frac{p-1}{p})}+t^2  \eta^{-\beta(\frac{p-1}{p})} \A^{\beta(\frac{p-1}{p})} +t^2 \sigma^{p-1} \bigg{\}}.
    \end{aligned}
\end{equation*} \noindent Now we divide left and right-hand side of this inequality for $t^{s+1}\rho^{N-s}$, we take the supremum on times $0<t<2\theta$ on the right and we integrate between $0\leq \tau \leq t$, to get
\begin{equation*}\label{passo11}
\begin{aligned}
%  \int_0^t \bigg[ \rho^{s-N}&\frac{d}{d\tau} \bigg(\tau^{-s}  \int_{B_t(y')}  \bigg( \int_{Q_{\tau,r}} u \zeta^p \, dx \bigg) \bigg]d\tau\leq\\
%  
\rho^{s-N} \bigg[t^{-s}&  \int \int_{\Q_{t,r}} u \zeta \, dx-\lim_{\tau \downarrow 0}  \dashint_{B_t(y')}  \bigg( \int_{B_r(\bar{x}'')} u \zeta dx''\bigg)dx'\bigg]\leq\\
%  &= \textcolor{violet}{\rho^{s-N} \bigg[ t^{-s} \int \int_{\Q_{t,r}} u \zeta \, dx-\int_{B_r(\bar{x}'')} u(y',\cdot) \zeta dx''\bigg]\leq} \quad \text{use} \, |\Q_{t,r}|= \gamma t^s r^{N-s}\approx  \gamma t^2 \rho^{N-s} \\
%  &\int_0^t \bigg{\{} \sup_{0<t<2\theta} \bigg[ t^{-1}\textcolor{blue}{t^{-s} \rho^{s-N}} \bigg( \frac{\gamma}{\rho \sigma^p} \textcolor{blue}{|\Q_{t,r}|} \bigg{\{} t^{\frac{2}{p}} \A^{2(\frac{p-1}{p})}+t^2  \eta^{-\beta(\frac{p-1}{p})} \A^{\beta(\frac{p-1}{p})} +t^2 \sigma^{p-1} \bigg{\}} \bigg)   \bigg] \bigg{\}}dt=\\
%  &\text{there remains a t downstairs but the total power of t is positive when we pass to the sup}\\
%  =
& \frac{\gamma t}{\rho \sigma^p} \bigg{\{} \theta^{\frac{2}{p}-1} \sup_{0<t<2\theta}\A^{2(\frac{p-1}{p})}+\theta  \eta^{-\beta(\frac{p-1}{p})} \sup_{0<t<2\theta}\A^{\beta(\frac{p-1}{p})} +\theta \sigma^{p-1} \bigg{\}}.
%  &\frac{\gamma}{\rho \sigma^p} \bigg{\{} \theta^{\frac{2}{p}} \sup_{0<t<2\theta}\A^{2(\frac{p-1}{p})}+\theta^2 \eta^{-\beta(\frac{p-1}{p})} \sup_{0<t<2\theta}\A^{\beta(\frac{p-1}{p})} +\theta^2 \sigma^{p-1} \bigg{\}}.
\end{aligned}
\end{equation*}\noindent Finally we use that $y'$ is a Lebesgue point and that $t < 2 \theta$ to have the estimate
\begin{equation} \label{passo12}
\begin{aligned}
\dashint & \dashint_{\Q_{t,r}} u \zeta dx \leq  \frac{1}{|B_r|} \int_{B_r(\bar{x}'')} u(y',\cdot) dx''+\\
+& \frac{\gamma}{\rho \sigma^p} \bigg{\{} \theta^{\frac{2}{p}} \sup_{0<t<2\theta}\A^{2(\frac{p-1}{p})}+\theta^2 \eta^{-\beta(\frac{p-1}{p})} \sup_{0<t<2\theta}\A^{\beta(\frac{p-1}{p})} +\theta^2 \sigma^{p-1} \bigg{\}}.
\end{aligned}
\end{equation} \noindent

\vskip0.4cm \noindent 
{\small{\it STEP 2. INTEGRAL INEQUALITY \eqref{passo12} IMPLIES A NONLINEAR ITERATION}}
\vskip0.2cm \noindent 
We consider $\varepsilon \in (0,1)$ and use Young's inequality to \eqref{passo12} to get for the first term
\begin{equation*}
    \begin{aligned}
\frac{\gamma \theta^{\frac{2}{p}}}{\rho \sigma^p} \bigg[\bigg(\sup_{0<t<2\theta}& \dashint \dashint_{\Q_{t,r}} u \, dx \bigg)^{\frac{2(p-1)}{p}}   +\eta^{\frac{2(p-1)}{p}}   \bigg] \leq \\
&\varepsilon \sup_{0<t<2\theta} \dashint \dashint_{\Q_{t,r}} u \, dx + \gamma(\varepsilon) \bigg(\frac{\gamma \theta^{\frac{2}{p}}}{\rho \sigma^p} \bigg)^{\frac{p}{2-p}}+\frac{\gamma \theta^{\frac{2}{p}}}{\rho \sigma^p} \eta^{\frac{2(p-1)}{p}}.
\end{aligned}
\end{equation*} \noindent Then we require a condition $\beta< \frac{p}{(p-1)}$ to get for the second term
\begin{equation*}
    \begin{aligned}
\frac{\gamma \theta^2 }{\rho \sigma^p}\eta^{-\beta(\frac{p-1}{p})}& \bigg[ \sup_{0<t<2\theta}\bigg( \dashint \dashint_{\Q_{t,r}} u \, dx \bigg)^{\beta(\frac{p-1}{p})} + \eta^{\beta(\frac{p-1}{p})}\bigg] \leq\\
& \varepsilon \sup_{0<t<2\theta} \dashint \dashint_{\Q_{t,r}} u \, dx + \gamma(\varepsilon) \bigg( \frac{\theta^2 \eta^{-\beta (\frac{p-1}{p})}}{\rho \sigma^p} \bigg)^{\frac{p}{p-\beta(p-1)}}+ \frac{ \gamma \theta^2}{\rho \sigma^p}.
\end{aligned} \end{equation*} \noindent 
Finally with these specifications, formula \eqref{passo12} is majorized by 
\begin{equation}\label{passo13}
\begin{aligned}
\dashint  \dashint_{\Q_{t,(1-\sigma)r}} u  \, dx &\leq
% &\rho^{s-N} \int_{B_r(\bar{x}'')} u(y',\cdot) dx''+ \epsilon \sup_{0<t<2\theta} |\Q|^{-1} \int \int_{\Q_{t,r}} u \, dx +\\
% +& C(\epsilon) \bigg(\frac{\gamma \theta^{\frac{2}{p}}}{\rho \sigma^p} \bigg)^{\frac{p}{2-p}}+\frac{\gamma \theta^{\frac{2}{p}}}{\rho \sigma^p} \eta^{\frac{2(p-1)}{p}}+ C(\epsilon) \bigg( \frac{\theta^2 \eta^{-\beta (\frac{p-1}{p})}}{\rho \sigma^p} \bigg)^{\frac{p}{p-\beta(p-1)}}+ \frac{ \gamma \theta^2}{\rho \sigma^p}+\frac{\gamma \theta^2}{\rho \sigma}\leq\\
\rho^{s-N} \int_{B_{2\rho}(\bar{x}'')} u(y',\cdot) dx''+ \varepsilon \sup_{0<t<2\theta} \dashint \dashint_{\Q_{t,r}} u \, dx +\\
% +& C(\epsilon) \gamma \sigma^{-p^2/(2-p)} \eta+\gamma   \sigma^{-p} \eta+ C(\epsilon) \bigg( \frac{\theta^2 \eta^{-\beta (\frac{p-1}{p})}}{\rho \sigma^p} \bigg)^{\frac{p}{p-\beta(p-1)}}+ \frac{ \gamma \theta^2}{\rho \sigma^p}\leq\\
+&\gamma \varepsilon^{-\tilde{\gamma}} \sigma^{-\tilde{\gamma}} \eta \bigg{\{} 1+ \eta^{-1}\bigg( \frac{\theta^2 \eta^{-\beta (\frac{p-1}{p})}}{\rho } \bigg)^{\frac{p}{p-\beta(p-1)}} + \eta^{-1} \frac{\theta^2}{\rho} \bigg{\}},
  \end{aligned}
\end{equation} \noindent being $\tilde{\gamma}>0$ a constant depending only on the data. The term in parenthesis is smaller than one if we contradict condition \eqref{either}, and the right-hand side of \eqref{passo13} is estimated from above with
\begin{equation}\label{passo14}
\begin{aligned}
 \dashint  \dashint_{\Q_{t,(1-\sigma)r}} u  dx \leq\rho^{s-N} \int_{B_r(\bar{x}'')} u(y',\cdot) dx''+ \varepsilon \sup_{0<t<2\theta} \dashint \dashint_{\Q_{t,r}} u \, dx + \gamma \varepsilon^{-\tilde{\gamma}} \sigma^{-\tilde{\gamma}} \eta.
  \end{aligned}
\end{equation} \noindent Since $t \in (0,2\theta)$ is an arbitrary number we choose the time $\bar{t}$ that achieves the supremum on the left of previous formula, 
% \[
% \sup_{0<t<2\theta}|\Q|^{-1}\int  \int_{\Q_{t,(1-\sigma)r}} u  dx=|\Q_{\bar{t},(1-\sigma)r}|^{-1}\int  \int_{\Q_{\bar{t},(1-\sigma)r}} u  dx,
% \]
and we observe that the right hand side of \eqref{passo14} does not depend on $t$. This implies 
\begin{equation*}\label{passo15}
\begin{aligned}
\sup_{0<t<2\theta} \dashint    \dashint_{\Q_{t,(1-\sigma)r}} u  dx \leq \rho^{s-N} \int_{B_{2\rho}(\bar{x}'')} u(y',\cdot) dx''+ \varepsilon \sup_{0<t<2\theta}\dashint \dashint_{\Q_{t,r}} u \, dx + \gamma \varepsilon^{-\tilde{\gamma}} \sigma^{-\tilde{\gamma}} \eta.
  \end{aligned}
\end{equation*} \noindent
We consider the increasing sequence $\rho_n=(1-\sigma_n)r= \rho (\sum_{i=1}^n 2^{-i}) \rightarrow 2 \rho$ and the expanding polydisc $Q_n=\{ Q_{t,\rho_n}\}: Q_{t,\rho} \rightarrow Q_{t,2\rho}$. By generality of the choice of $r$, previous formula implies the recurrence 
\begin{equation}\label{passo16}
\S_n:= \sup_{0<t<2\theta} \dashint  \dashint_{\Q_n} u  dx \leq
\varepsilon \S_{n+1} + \bigg( \rho^{s-N} \int_{B_{2\rho}(\bar{x}'')} u(y',\cdot) dx'' \gamma \varepsilon^{-\tilde{\gamma}} \eta \bigg) b^{n}, 
\end{equation} \noindent with $b= 2^{\tilde{\gamma}}$. Therefore with an iteration as in \eqref{sqn} we arrive to the conclusion
\begin{equation}\label{passo17}
\begin{aligned}
\sup_{0<t<2\theta} \dashint  \dashint_{\Q_{t,\rho}} u  \,dx \leq \rho^{s-N} \int_{B_{2\rho}(\bar{x}'')} u(y',\cdot) dx''+  \gamma \eta.
  \end{aligned}
\end{equation} \noindent \vskip0.2cm \noindent 
{\small{\it STEP 3. ESTIMATING THE FULL $L^{\infty}$} NORM FROM ABOVE.} \vskip0.1cm 
\noindent We consider formula \eqref{Ls-Linfty} with $l=1$, which is
\begin{equation}\label{MYTHOLOGY}
\sup_{Q_{\theta/2},\rho/2} u 
 \leq \gamma  \bigg( \frac{\rho^p}{\theta^2} \bigg)^{\frac{(N-s)}{\chi}}  \bigg( \dashint  \dashint_{Q_{\theta,\rho}} u_+\, dx  \bigg)^{\frac{p}{\chi}} +\gamma \bigg(\frac{\theta^2}{\rho^p} \bigg)^{\frac{1}{2-p}}.
\end{equation} \noindent 
Let us insert in the integral term on the right of \eqref{MYTHOLOGY} our previously obtained formula \eqref{passo17} to get
\begin{equation}\label{titto}
\begin{aligned}
\sup_{Q_{\theta/2},\rho/2} u &
 \leq \gamma  \bigg( \frac{\rho^p}{\theta^2} \bigg)^{\frac{(N-s)}{\chi}}  \bigg( \rho^{s-N} \int_{B_{2\rho}(\bar{x}'')} u(y',\cdot)+\gamma \eta \bigg)^{\frac{p}{\chi}} +\gamma \bigg(\frac{\theta^2}{\rho^p} \bigg)^{\frac{1}{2-p}}\\
 &\leq \gamma  \bigg( \frac{\rho^p}{\theta^2} \bigg)^{\frac{(N-s)}{\chi}}  \bigg( \rho^{s-N} \int_{B_{2\rho}(\bar{x}'')} u(y',\cdot) \bigg)^{\frac{p}{\chi}} +\gamma \bigg(\frac{\theta^2}{\rho^p} \bigg)^{\frac{1}{2-p}}.
 \end{aligned}
\end{equation} \noindent The proof is complete.
\vskip0.4cm \noindent

{\centering{
\section{Harnack inequality. Proof of Theorem \ref{harnackTHM}\label{7}}}}
\noindent 
Here we prove the Harnack estimate \eqref{Harnack}, by use of Theorems \ref{shrinkingTHM}, \ref{l1-linftyTHM}. Without loss of generality we assume that $x_0=0$ and denote $u(x_0)=u_0$ to ease notation. First we begin with an estimate reminescent of \cite{KS} aiming to a bound from above and below in terms of the radius itself (estimate \eqref{sup-estimate}).  \noindent \vskip0.2cm \noindent 
{\small{\it STEP 1. A Krylov-Safonov argument}.} \vskip0.1cm

\begin{figure}
\begin{tikzpicture}[scale=0.4]
\draw (0,0) circle (6cm) node{$O$};
\draw[dashed] (-0.4,0) -- (-6,0) node[above right]{$\lambda_1$};
\draw (0,0) circle (3cm);
\draw[dashed] (0,-0.4) -- (0,-3)node[below right]{$\lambda_0$};
\draw[dashed] (1.5,1.8) -- (1.5,4) node[below right]{$2r$};
\draw (1.5,1.5) circle (2.5cm) node{$\bar{x}$};

\draw[thick,->] (15,-6) -- (27,-6) node[anchor=north west] {$\lambda$};
\draw[thick,->] (15,-6) -- (15,6);
\draw[black,thick,dashed] (25,-6) -- (25,6);
\draw (25,-6) node[below]{$1$};
\draw (14.5,0) node[below]{$1$};
\draw (15,0) node[below]{$-$};
\draw (26,0) .. controls (22,-6) and (18,10) .. (15,-0.5);

\draw[blue] (15,-0.5) .. controls (18,10) and (22,-6) .. (26,0);
\draw[blue] (16,3) node[above right]{$g(\lambda)$};
\draw[red] (15,-0.5) parabola (24,6);

\draw[red] (24,5) node[above left]{$f(\lambda)$};
\draw[black,thick,dashed] (20,1.6) -- (20,-6);
\draw (20,-6) node[below]{$\lambda_0$};
\end{tikzpicture}

\caption{\small{Figure on the left illustrates the geometric construction of $B_{2r}(\bar{x})$, with $\rho=1$ and $u_0=1$. On the right we draw the idea of definition of $\lambda_0$ through equation \eqref{Krylov-Safonov-EQ}. Function $f$ (in red) represents the member on the right-hand side of equation \eqref{Krylov-Safonov-EQ} while function $g$ (in blue) represents the left-hand side member.}}\label{figA}
\end{figure}
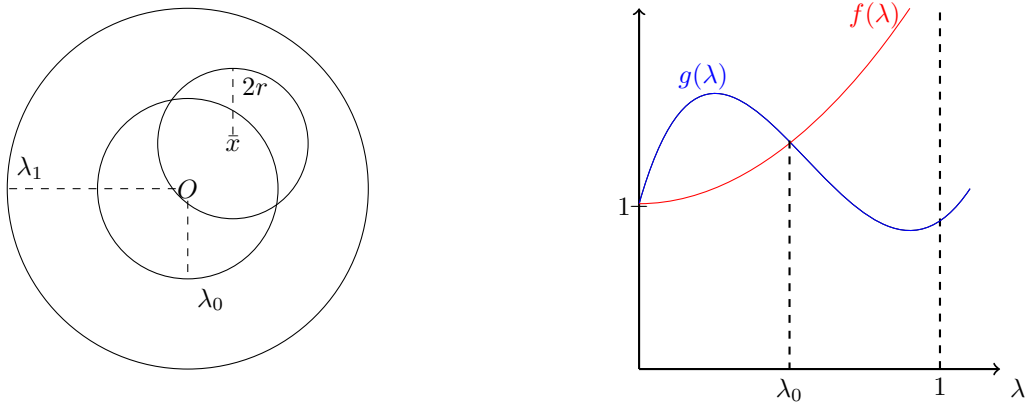

\noindent For a parameter $\lambda \in (0,1)$ we consider the equation
\begin{equation} \label{Krylov-Safonov-EQ}
    \sup_{x''\in B_{\lambda}\rho(0'')} u(0',x'') = u_0 (1-\lambda)^{-\beta}.
\end{equation} \noindent
Let $\lambda_0$ be the maximal root of the equation \eqref{Krylov-Safonov-EQ} and by continuity let us fix a point $\bar{x}''$ by
\begin{equation*} \label{thepoint}
    u(0',\bar{x}'')= \max_{x''\in B_{\lambda_0}\rho(0'')} u(0',x'') = u_0 (1-\lambda_0)^{-\beta}=:M, \end{equation*} being $\bar{x}''$ a point in $B_{\lambda_0 \rho}(0'')$.  See Figure \ref{figA} for a drawing.  Now let us define $\lambda_1 \in (0,1)$ by 

\noindent 

\noindent \[ (1-\lambda_1)^{-\beta}=4(1-\lambda_0)^{-\beta}, \quad \quad \text{i.e.} \quad  \quad \lambda_1=1-4^{-1/\beta}(1-\lambda_0), \quad \lambda_1 >\lambda_0,\] and set also 
\[
2r:= (\lambda_1-\lambda_0)\rho=(1-4^{-1/\beta})(1-\lambda_0)\rho.
\]Then by definition of $\lambda_0$ it holds both
\[
\sup_{x''\in B_{\lambda_1 \rho}(0'')} u(0', x'') \leq u_0(1-\lambda_1)^{-\beta}\quad  \quad \text{\&} \quad \quad  B_{2r}(\bar{x}'') \subset B_{\lambda_1\rho}(0'').
\]

\noindent See Figure \ref{figA} for this construction. Henceforth we arrive to the estimate
\begin{equation} \label{sup-estimate}
\begin{aligned}
    % \sup_{x'' \in B_{\lambda_0 \rho}(\bar{x}'')} u(0',x'') 
    % \frac{M}{2}=&(1-\lambda_0)^{-\beta}=u(0',\bar{x}'') \leq 
    M=u(0',\bar{x}'') \leq \sup_{x'' \in B_{2r}(\bar{x}'')}& u(0',x'') \leq
    % \\
    % & \sup_{x''\in B_{\lambda_1 \rho}(0'')} u(0',x'') \leq
    u_0 (1-\lambda_1)^{-\beta}= 4 u_0(1-\lambda_0)^{-\beta}=4M.
    \end{aligned}
\end{equation}  

\noindent Now we use the bound obtained joint to the $L^1-L^{\infty}$ estimate \eqref{l1-linfty} to reach an estimate of the measure of sub-level sets of $u$, in order to apply Theorem \ref{shrinkingTHM}. 
\noindent \vskip0.2cm \noindent 
{\small{\it STEP 2. Estimating the slice-wise measure of sub-levels of $u$}.} \vskip0.1cm

\noindent Further we assume $u_0 \ge K \rho$, where $K$ is the number of Theorem \ref{shrinkingTHM}. Let us construct the polydisc 
\[
Q_{\eta,2r}(0',\bar{x}''), \quad \quad \text{with} \quad \quad  \eta= (2r)^{\frac{p}{2}}M^{\frac{2-p}{2}}, \quad \quad  \text{being} \quad \quad  Q_{\eta,2r}(0,\bar{x}'') \subset Q_{\mathcal{M},\rho}(0)\subset \Omega.
\] \noindent Last inclusions are due to the fact that $r<\rho$, $M \leq 2 ||u||_{L^{\infty}(\Omega)}$ and hypothesis \eqref{M}. Hence we apply Theorem \ref{l1-linftyTHM} in $Q_{\eta,2r}(0,\bar{x}'')$ and use \eqref{sup-estimate}  to get the inequality
\begin{equation}\label{l1}
\begin{aligned}
\max_{Q_{\frac{\eta}{2},\frac{r}{2}}(0',\bar{x}'')}u \leq& \gamma \bigg{\{} \bigg( \frac{r^p}{\eta^2} \bigg)^{\frac{N-s}{\chi}} \bigg[ \inf_{x'\in B_{\eta}(0')}  \dashint_{B_{2r}(\bar{x}'')} u(x',x'')\, dx''  \bigg]^{\frac{p}{\chi}}  + \bigg( \frac{\eta^2}{r^p} \bigg)^{\frac{1}{2-p}} \bigg{\}}\leq \gamma M.
\end{aligned}
\end{equation} \noindent Next, we use again Theorem \ref{l1-linftyTHM} for the polydisc $Q_{\delta \eta,r/4 }(0,\bar{x}'')$ for a $\delta>0$ to be determined later. By choice of $\bar{x}''$ we obtain the estimate
\begin{equation}\label{l2}
\begin{aligned}
M=u(0',\bar{x}'')\leq \gamma \bigg{\{} \bigg(\frac{M^{p-2}}{\delta^{2}}\bigg)^{\frac{(N-s)}{\chi}} \bigg[ \inf_{x'\in B_{\delta \eta}(0')}  \dashint_{B_{r/2}(\bar{x}'')} u(x',x'')\, dx''  \bigg]^{\frac{p}{\chi}}  + \delta^{\frac{2}{2-p}}M \bigg{\}},
\end{aligned}
\end{equation}\noindent and we estimate the integral term by splitting the integral on the level $\tilde{\epsilon}$ as
\begin{equation*}
    \begin{aligned}  \inf_{x' \in B_{\delta \eta}(0')} r^{s-N} \bigg(  \int_{B_{r/2}(\bar{x}'') \cap [u \ge \tilde{\epsilon} M]}u(x',x'') \, dx''+ \tilde{\epsilon} M | B_{r/2}(\bar{x}'')\cap[u(x',\cdot) <\tilde{\epsilon} M]|  \bigg).
    \end{aligned}
\end{equation*}\noindent First integral term is estimated by \eqref{l1} and making the choice $\gamma \delta^{\frac{2}{2-p}}= 1/4$ we have from \eqref{l2} that
\begin{equation}\label{l3}
    \begin{aligned}
    \frac{M}{2}\leq \gamma \delta^{\frac{2(s-N)}{\chi}}  M \bigg(\tilde{\epsilon} +  
    \frac{\inf_{x'\in B_{\delta \eta}(0')} |[u(x',\cdot) \ge \tilde{\epsilon} M] \cap B_{r/2}(\bar{x}'')|}{|B_{r/2}(\bar{x}'')|}   \bigg)^{\frac{p}{\chi}}+ M/4.
        \end{aligned}
\end{equation} \noindent  Henceforth, by choosing also $\tilde{\epsilon} \leq  \frac{1}{2}\bigg( \frac{\delta^{2(N-s)}}{(4\gamma)^{\chi}}\bigg)^{1/p}$ we obtain the measure estimate
\[
 \frac{|u(x',\, \cdot) \ge \tilde{\epsilon} M] \cap B_{r/2}(\bar{x}'')|}{|B_{r/2}(\bar{x}'')|} \ge  \frac{1}{2}\bigg( \frac{\delta^{2(N-s)}}{(4\gamma)^{\chi}}\bigg)^{1/p}=:\alpha,
\] \noindent for each $x' \in B_{\delta \eta}(0')$. This gives us the inequality
\begin{equation}\label{alzacheschiaccio}
|[u(x',\cdot) \leq \tilde{\epsilon} M] \cap B_{r/2}(\bar{x}'')| \leq (1-\alpha) |B_{r/2}(\bar{x}'')|, \quad \forall_{ae}\, x' \in B_{\delta \eta}(0').
\end{equation} \noindent If we reduce furthermore $\tilde{\epsilon} < \min\{\tilde{\epsilon}, (\delta 2^{p/2})^{\frac{2}{2-p}}  \} $ and use that $1<p<2$, then we have
\[\theta:= (r/2)^{\frac{p}{2}} (\tilde{\epsilon} M)^{\frac{2-p}{2}}\leq \delta \eta = \delta r^{\frac{p}{2}} M^{\frac{2-p}{2}},\] just in order to apply Theorem \ref{shrinkingTHM} in the polydisc $Q_{\theta,r/2}(0',\bar{x}'')$. \noindent Finally here below we expand the positivity applying iteratively Theorem \ref{shrinkingTHM} to $u$ in appropriate neighborhods of $(0',\bar{x}'')$, in order to expand positivity until we reach a neighborhood of the origin. A lower bound which is free from any dependence on $u$ itself can be achieved by choosing $\beta$ appropriately.
\noindent \vskip0.2cm \noindent 
{\small{\it STEP 3. Expansion of positivity and choice of $\beta$}.} \vskip0.1cm 

\noindent We consider the measure estimate \eqref{alzacheschiaccio}: either holds \eqref{alternative-harnack} or 
\[
u(x) >\delta_o \tilde{\epsilon} M/2, \quad \text{in} \quad  Q_{\eta,r}(0',\bar{x}''), \quad \text{being} \quad \eta=(r/2)^{\frac{p}{2}}(\delta_o \tilde{\epsilon} M)^{\frac{2-p}{2}}.\]
This implies the measure estimate
\begin{equation}\label{m1}
    |[u(x',\cdot)\leq \delta_o \tilde{\epsilon}M/2] \cap B_{r}(\bar{x}'')|\leq |B_{r}|/2.
\end{equation}\noindent 
We can apply again Theorem \ref{shrinkingTHM}. This time and next ones being $\alpha=1/2$ fixed, there exists a number $\delta_*$ depending only on the data and such that for almost every $x \in Q_{\eta_*,2r}(0',\bar{x}'')$ we have
\begin{equation}\label{m2}
    u(x)> \delta_* M^*/2, \quad\quad\eta_*^2= (2r)^{p}(\delta_* M^*)^{2-p}  \quad  \quad \text{being} \quad M^*= \delta_o \tilde{\epsilon}M/2.
\end{equation}\noindent Now the procedure can be iterated a number $n \in \N$ of times such that $n \ge\log_2(4/(1-\lambda_0))$ in order to have 
\[
u(x)> (\delta_*/2)^{n}M^{*}, \quad \forall_{ae} \, x \in Q_{\eta_*(n),2^n r}(0',\bar{x}''), \quad \text{being}\quad \eta_*^2(n)=(2r)^p(\delta_*^n M^*)^{2-p}.
\] We observe that in previous calculation of $\eta_*(n)$ the powers of $2$ cancel each other out.\newline
\noindent Since $\bar{x}'' \in B_{\lambda_0 \rho}(0'')$, then $B_{\rho}(0'') \subset B_{2\rho}(\bar{x}'')$. If we assume $\beta>2$ the choice of $\bar{n}$ above implies $2^{\bar{n}} r \ge 2\rho$ so that 
\begin{equation} \label{fine}
u(x) > (\delta_*/2)^{n}M^*, \quad \quad \forall_{ae}\, x \in Q_{\eta_*(n), \rho}(0), \quad   \text{where}
\end{equation}
\begin{equation*} \begin{aligned}
\eta_*(n)=&(2 r)^{\frac{p}{2}}(\delta_*^n M^*)^{\frac{2-p}{2}}=\\
&
(2r)^{\frac{p}{2}}(\delta_*^n \delta_o \tilde{\epsilon} M)^{\frac{2-p}{2}}= (1-\lambda_0)^{p/2}\rho^{p/2}\{ \delta_*^n\delta_o \tilde{\epsilon} [u_0(1-\lambda_o)^{-\beta}2^{-1}]\}^{\frac{2-p}{2}}\ge \bar{\delta}_o \rho^{p/2} u_0^{\frac{2-p}{2}},
\end{aligned} \end{equation*}\noindent  by properly choosing $\beta\ge p/(2-p)$ and redefining the constants. Observe that equation \eqref{fine} is exactly \eqref{Harnack}. To end the proof, we will choose $\beta$ big enough to free the lower bound 
\[u(x) > (\delta_*/2)^{n}M^*,\] by any dependence of the solution itself other than $u_0$. Indeed, decreasing $\rho$ in case of need, let $n\in \N$ be a number big enough that
\[
n \ge \bar{n}, \quad \quad 1\leq 2^n \bigg( \frac{r}{\rho}\bigg) \leq 2, \quad \Rightarrow \quad (1-\lambda_0)^{-1} >2^{n-2}(1-4^{-1/\beta}).
\] Then we have 
 \begin{equation*} 
(\delta_*/2)^{n}M^*= (\delta_*/2)^{n}[\delta_0 \tilde{\epsilon} u_0 (1-\lambda_0)^{-\beta}]\ge (\delta_* 2^{\beta-1})^n 2^{-2\beta} (1-4^{-1/\beta})^{\beta} \delta_0 \tilde{\epsilon} \, u_0. \end{equation*} 
\noindent  Decreasing $\delta_*$ in case of need, we choose finally $\beta>p/(2-p)$ so big that 
\[\delta_* 2^{\beta-1}=1, \quad \quad K:=2^{-2\beta} (1-4^{-1/\beta})^{\beta} \delta_0 \tilde{\epsilon} \] and the claim follows.

% \newpage \noindent 
% \begin{remark}
% \textcolor{blue}{Let us consider the measure information \eqref{alzacheschiaccio}, which is 
% \begin{equation}\label{alzacheschiaccio}
% |[u(x',\cdot) \leq \tilde{\epsilon} M] \cap B_{r/2}(\bar{x}'')| \leq (1-\alpha) |B_{r/2}(\bar{x}'')|, \quad \forall_{ae}\, x' \in B_{\delta \eta}(0').
% \end{equation} \noindent This implies directly
% \begin{equation*}
% |[u(x',\cdot)\leq \tilde{\epsilon}M \cap B_{2^{\bar{n}+1}r}(\bar{x}'')]|\leq \bigg(1-\frac{\alpha}{2^{\bar{n}(N-s)}}  \bigg)    |B_{2^{\bar{n}+1}r}|,
% \end{equation*} \noindent for each $x'\in B_{\delta \eta} (0')$ and for $\bar{n}$ defined as in the proof above. Now we just use just use once Theorem \ref{shrinkingTHM} to get the existence of $K>1,\, \delta_o\in (0,1)$ depending only on the data and $\alpha$ such that
% \begin{equation*}
% u(x) > \delta_o\tilde{\epsilon}M/2, \quad \forall_{ae} \, x \in Q_{\eta(n),2^{\bar{n}+2}r}(0',\bar{x}'').
% \end{equation*}}
% \end{remark}

% Here the idea is that it is not important in Theorem \ref{shrinkingTHM} that we expand positivity. We can just do it now, when we finish the proof of the Harnack.

% \newpage

{\centering{\section{Appendix}\label{8}}}
\noindent 
In this Appendix we enclose all the details that for reader's convenience have been postponed. \vskip0.1cm \noindent 

\subsection{Proof of the Energy Estimates \eqref{energy1}} \label{e1proof}
We suppose without loos of generality that $\bar{x}=0$. Let us test the equation \eqref{def-solution} with $\phi= \pm (u-k)_{\pm} \zeta^2 \in W_o^{1,[2,p]}(\Omega)$ for $\zeta \in C_o^{\infty}(\Omega)$ to obtain
\begin{equation} \label{ee-a}
    \begin{aligned}
    0=\int \int_{Q_{\theta,\rho}} \bigg{\{}& \sum_{i=1}^s (\partial_i u) \bigg[\bigg(\pm \partial_i (u-k)_{\pm}\bigg) \zeta^2 +2(u-k)_{\pm} (\partial_i \zeta) \zeta      \bigg]+\\
    &+\sum_{i=s+1}^N A_i(x,\nabla u) \bigg[\bigg(\pm \partial_i (u-k)_{\pm}\bigg) \zeta^2 +2(u-k)_{\pm} (\partial_i \zeta) \zeta      \bigg] \bigg{\}} dx.
    \end{aligned}
\end{equation}\noindent We divide the terms in squared parenthesis and use \eqref{E1-structure} to get
\begin{equation}\label{ee-b}
    \begin{aligned}
    \sum_{i=1}^s &\int \int_{Q_{\theta,\rho}} |\partial_i (u-k)_{\pm}|^2 \, dx + C_1 \sum_{i=s+1}^N \int \int_{Q_{\theta,\rho}} |\partial_i (u-k)_{\pm}|^p \, dx \leq \\
    & (N-s) C \int \int_{Q_{\theta,\rho}} \chi_{[u\lesseqgtr k]}(x) dx+  \sum_{i=1}^s \int \int_{Q_{\theta,\rho}} \epsilon |\partial_i u|^2 \zeta^2 + C(\epsilon) |(u-k)_{\pm}|^2 |\partial_i \zeta|^2 dx \\
    &+2 \sum_{i=s+1}^N \int \int_{Q_{\theta,\rho}} \bigg[ C_2 |\partial_i (u-k)_{\pm}|^{p-1}+C \bigg] \bigg[ (u-k)_{\pm} (2\zeta)\partial_i \zeta \bigg]\, dx,
    \end{aligned}
\end{equation} \noindent where we have used Young inequality $ab \leq \epsilon a^2+C(\epsilon) b^2$ on the term $ (\partial_i u) [ 2(u-k)_{\pm} (\partial_i \zeta) \zeta]$ of the first line of \eqref{ee-a}.\newline 
Two other similar applications of Young's inequality $ab \leq \epsilon a^p+C(\epsilon) b^{\frac{p}{p-1}}$ to the last term on the right of \eqref{ee-b} reveal that this is smaller than the quantity
\[
2 \sum_{i=s+1}^N \int \int_{Q_{\theta,\rho}} \bigg[\epsilon |\partial_i (u-k)_{\pm}|^p+2 C(\epsilon) |(u-k)_{\pm}|^p \, |\partial_i \zeta|^p + \epsilon (2 \zeta)^p \chi_{\{(u-k)_{\pm}\ge 0\}}(x)  \bigg] \, dx.
\] Gathering together all this estimates and choosing $\epsilon>0$ to be a constant small enough concludes the proof.

\subsection{Iteration Lemmata}
Here we recall two basic Lemmas, extremely useful for the iteration techniques employed in our analysis and whose proofs can be found in (\cite{DB}, Chap. I Sec.IV).
\begin{lemma}[\cite{DB} Chap. I, Sec. IV] \label{Fast}
Let $\{ Y_n\}_{n \in \N}$ be a sequence of positive numbers satisfying the recursive inequalities
\begin{equation}
    Y_{n+1} \leq C b^n Y_{n}^{1+\alpha},
\end{equation} \noindent where $C,b>1$ and $\alpha \in (0,1)$ are given numbers. Then we have the logical implication
\[
Y_o \leq C^{\frac{-1}{\alpha}} b^{\frac{-1}{\alpha^2}} \quad \Rightarrow \quad \lim_{n \rightarrow \infty} Y_n=0.
\]
\end{lemma} \noindent
The following Lemma is useful for reverse recursive inequalities of Section \ref{DG-IT-proof}.

\begin{lemma}[\cite{DB} Chap. I, Sec. IV]\label{reverse}
Let $\{ Y_n \}_{n \in \N}$ be a sequence of equibounded positive numbers which satisfies the recursive inequality
\begin{equation}\label{recursive-reverse}
    Y_n \leq C b^n Y_{n+1}^{1-\alpha}, \quad \quad \alpha \in (0,1), \quad C,\, b>1.
\end{equation} \noindent Then
\begin{equation}\label{TH-reverse}
    Y_0 \leq \bigg( \frac{2C}{b^{1-1/\alpha}}  \bigg)^{\frac{1}{\alpha}}.
\end{equation}
\end{lemma}

% \begin{proof}
% Let us use Young's inequality to \eqref{recursive-reverse} to get
% \begin{equation} \label{point-of-not-return}
% Y_n\leq \epsilon Y_{n+1}+ \bigg(\frac{C}{\epsilon^{1-\alpha}}  \bigg)^{\frac{1}{\alpha}} b^{\frac{n}{\alpha}}, \quad \forall \epsilon \in (0,1), \quad n=0,1,2,...
% \end{equation} \noindent Let us do one step to understand the iteration
% \[
% Y_{n-1}\leq \epsilon Y_n+\bigg(\frac{C}{\epsilon^{1-\alpha}}  \bigg)^{\frac{1}{\alpha}} b^{\frac{n-1}{\alpha}}\leq  \epsilon \{ \epsilon Y_{n+1}+\bigg(\frac{C}{\epsilon^{1-\alpha}}  \bigg)^{\frac{1}{\alpha}} b^{\frac{n}{\alpha}} \}+\bigg(\frac{C}{\epsilon^{1-\alpha}}  \bigg)^{\frac{1}{\alpha}} b^{\frac{n-1}{\alpha}}
% \] so that by iteration we get
% \[
% Y_o \leq \epsilon^n Y_n + \bigg( \frac{C}{\epsilon^{1-\alpha}}  \bigg)^{\frac{1}{\alpha}} \sum_{i=0}^{n-1} (b^{\frac{1}{\alpha}} \epsilon)^{i}.
% \] Finally choosing $b^{\frac{1}{\alpha}}\epsilon=\frac{1}{2}$ the sum on the right side can be majorised by the corresponding series converging to 2. Letting $n \rightarrow \infty$ finishes the prof.
% \end{proof}

\begin{remark} \label{rmk}
If we just have a sequence of equibounded positive numbers $\{Y_n\}$ such that
\begin{equation}\label{sqn}
    Y_n \leq \epsilon Y_{n+1}+ C b^n, \quad C,b>1, \quad \epsilon \in (0,1),
\end{equation} \noindent then by a simple iteration, setting $(\epsilon b)=1/2$ and letting $n \rightarrow \infty$ gives \eqref{TH-reverse} with $\alpha=1$. See for instance (\cite{DB}, Lemma 4.3 page 13).
\end{remark}

\subsection{Proof of Lemma \ref{DG-IT}} \label{DG-IT-proof}
\begin{proof} We are going to perform a cross-iteration. Let $\{\sigma_j\}_{j \in \N} \subset (0,1)$ be the increasing sequence $\sigma_j= \sum_{i=1}^j 2^{-(i+1)}$. Let also $n \in \N$ be an index to define the decreasing sets
\begin{equation*}
\begin{cases} Q_n= Q_{\theta_n,\rho_n}, \quad \tilde{Q}_n=Q_{\tilde{\rho}_n,\tilde{\theta}_n},\\
Q_0=Q_{\sigma_{j+1}\theta,\sigma_{j+1} \rho}, \quad Q_{\infty}= Q_{\sigma_j \theta, \sigma_j \rho},
\end{cases} \text{for} \quad 
\begin{cases}
\rho_n= \sigma_{j} \rho+ \frac{(\sigma_{j+1}-\sigma_j)\rho}{2^n}, \quad \tilde{\rho_n}= \frac{\rho_n+\rho_{n+1}}{2},\\
\theta_n= \sigma_j \theta+ \frac{(\sigma_{j+1}-\sigma_j)\theta}{2^n}, \quad \tilde{\theta}= \frac{\theta+\theta{n+1}}{2}.
\end{cases}
\end{equation*} \noindent
Furthermore, contradicting \eqref{raggae0}, let $k\ge ( \theta^2/\rho^p)^{\frac{1}{2-p}}\ge \rho$ be a number to be defined a posteriori and  let us define the increasing sequence of levels and numbers
\[
k_n=k-\frac{k}{2^n}, \quad \quad M_{n}= \sup_{Q_n} u,\]
such that
\[M_0=:M_{j+1}:= \sup_{Q_{\sigma_{j+1}\theta,\sigma_{j+1}\rho}} u, \quad \quad  M_{\infty}=:M_j:= \sup_{Q_{\sigma_j \theta, \sigma_j \rho}} u.
\]

 \noindent 
{\small{\it STEP 1. FIRST ITERATION (valid without condition \eqref{chi})}}
\vskip0.2cm 

\noindent We perform in the first place an iteration for $n \rightarrow \infty$ on shrinking polydiscs.
%  , until we will have an estimate of the kind
% \[
% E (Q_{\sigma_j \theta, \sigma_j \rho}) \leq E(Q_{\sigma_{j+1}\theta, \sigma_{j+1}\rho}),
% \] and then an iteration on $\sigma_j$ on expanding cylinders, such that $\sigma_0=1/2$ and $\sigma_n \rightarrow 1$, which will lead us to define a sequence $ \{Y_n\}$ and use Lemma \ref{reverse}.
% \end{remark}
\noindent To this aim, we introduce cut-off functions $\zeta_n$ vanishing on $\partial Q_n$ and equal to one in $\tilde{Q}_n$, that obey to 
\begin{equation*}
\begin{cases}
|\partial_i \zeta_n|\leq (\rho_{n}-\rho_{n+1})^{-1}= (2^{n+1} 2^{j+1})/\rho,& \forall i=1,..,s,\\
|\partial_i \zeta_n|\leq (\theta{n}-\theta{n+1})^{-1}=(2^{n+1} 2^{j+1})/\theta,&\forall i=s+1,..,N.
\end{cases}
\end{equation*}\noindent 
With these stipulations,the energy estimates \eqref{energy1} are
\begin{equation}\label{energade}
\begin{aligned}
&I_n:= \int \int_{\tilde{Q}_n} \sum_{i=1}^s | \partial_i (u-k_{n+1})_{+}|^{2} \zeta_n^2 \,  + \sum_{i=s+1}^N  | \partial_i (u-k_{n+1})_{+}|^{p} \zeta_n^2 \, dx\\
&\leq     \gamma \int \int_{Q_n} \bigg{\{}  \frac{ 2^{2n}2^{2j}}{ \theta^2}   |(u-k_{n+1})_{+}|^{2}  + \frac{ 2^{np}2^{jp}}{ \rho^p}  | (u-k_{n+1})_{+}|^{p}  \bigg{\}}\, dx+ |A_{n}|,
\end{aligned}
\end{equation}\noindent being $A_n= Q_n\cap [u> k_{n+1}]$. Now for any $s>0$ we observe that
\[
|A_n|\leq \frac{2^{s(n+1)}}{k^s}\int \int_{Q_n} (u-k_{n})_+^s dx.
\]
This fact with $s=2$ together with H\"older inequality gives
\begin{equation*}
    \begin{aligned}
\int \int_{Q_n} (u-k_{n+1})^p_+\,dx \leq& \bigg( \int \int_{Q_n} (u-k_{n+1})_+^2\, dx  \bigg)^{\frac{p}{2}} |A_n|^{1-\frac{p}{2}}\\
\leq & \gamma \frac{2^{(2-p)n}}{k^{2-p}} \int \int_{Q_n} (u-k_{n})_+^2 dx.
    \end{aligned}
\end{equation*} \noindent Putting this into the energy estimates above leads us to
\begin{equation}\label{raggae1}
    \begin{aligned}
    I_{n}:=\int \int_{\tilde{Q}_n}& \bigg{\{}\sum_{i=1}^s |\partial_i (u-k_{n+1})_+|^2 + \sum_{i=s+1}^N |\partial_i (u-k_{n+1})_+|^p\bigg{\}} dx\\
    \leq&\frac{\gamma 2^{2n}2^{2j}}{\theta^2} \bigg(1+k^{p-2}\frac{\theta^2}{\rho^p} +\frac{\theta^2}{k^2} \bigg) \int \int_{Q_n} |(u-k_n)_+|^2\, dx\\
    % \leq & \frac{\gamma 2^{2n}}{(1-\sigma)^2\theta^2} \bigg(1+k^{p-2}\frac{\theta^2}{\rho^p}\bigg) \int \int_{Q_n} |(u-k_n)_+|^2\, dx\\
    &\leq \frac{\gamma 2^{2n}2^{2j}}{\theta^2} \int \int_{Q_n} |(u-k_n)_+|^2\, dx,
    \end{aligned}
\end{equation} \noindent 
because $k\ge (\theta^2/\rho^p )^{\frac{1}{2-p}}\ge \rho$. Now an application of Troisi's Lemma \ref{Troisi} and \eqref{raggae1} above give us the following inequality

\begin{equation}
    \begin{aligned}
    \int \int_{Q_{n+1}}& (u-k_{n+1})_+^2\, dx  \leq \int \int_{\tilde{Q}_n}  (u-k_{n+1})_+^2 \tilde{\zeta}^2\, dx\leq \\
    &\leq \bigg( \int \int_{Q_n} [(u-k_{n+1})_+ \tilde{\zeta}]^{\frac{N\bar{p}}{N-\bar{p}}} \bigg)^{2(\frac{N-\bar{p}}{N\bar{p}})} |A_n|^{1-\frac{2}{\bar{p}}+\frac{2}{N}}\\
    &\leq \gamma \bigg( \prod_{i=1}^N ||\partial_i (u-k_{n+1})_+||_{L^{p_i}(Q_n)} \bigg)^{\frac{2}{N}}|A_n|^{1-\frac{2}{\bar{p}}+\frac{2}{N}} \\
    &\leq \gamma \bigg( I_n^{\sum_{i=1}^N \frac{1}{p_i}}  \bigg)^{\frac{2}{N}}|A_n|^{1-\frac{2}{\bar{p}}+\frac{2}{N}} \,\\
    &\leq \bigg[ \frac{\gamma 2^{2n}2^{2j}}{\theta^2} \bigg]^{\frac{2}{\bar{p}}} \bigg( \frac{2^{2(n+1)}}{k^2} \bigg)^{1-\frac{2}{\bar{p}}+\frac{2}{N}}  \bigg( \int \int_{Q_n} |(u-k_n)_+|^2\, dx\bigg)^{1+\frac{2}{N}}.
    \end{aligned}
\end{equation} \noindent 
Hence by setting
\[
Y_n= \int_{Q_n} (u-k_n)_+^2\, dx, 
\] we arrive at the inequality
\[
Y_{n+1} \leq \gamma 2^{n(1+\frac{2}{N})} 2^{\frac{4j}{\bar{p}}} \theta^{-\frac{4}{\bar{p}}} k^{-2(\frac{N(\bar{p}-2)+2\bar{p}}{N\bar{p}})} Y_{n}^{1+\frac{2}{N}},
\] which converges by Lemma \ref{Fast} if
\[
Y_0\leq \bigg( \gamma 2^{\frac{4j}{\bar{p}}}  \theta^{-\frac{4}{\bar{p}}}k^{-2(\frac{N(\bar{p}-2)+2\bar{p}}{N\bar{p}})} \bigg)^{-\frac{N}{2}} 2^{\frac{N+2}{N}(-\frac{N^2}{4})}= \gamma 2^{-\frac{2Nj}{\bar{p}}} \theta^{\frac{2N}{\bar{p}}} k^{\frac{N(\bar{p}-2)+2\bar{p}}{\bar{p}}}.
\]
This condition can be obtained by imposing
\begin{equation}\label{k-in-sup}
    k= \gamma 2^{\tilde{\gamma}j} \theta^{\frac{-2N}{N(\bar{p}-2)+2\bar{p}}} \bigg( \int \int_{Q_{\sigma_{j+1}\theta,\sigma_{j+1}\rho}} u_+^2\, dx  \bigg)^{\frac{\bar{p}}{N(\bar{p}-2)+2\bar{p}}} \wedge \bigg(\frac{\theta^2}{\rho^p} \bigg)^{\frac{1}{2-p}}.
\end{equation} \noindent Therefore we proceed by this choice of $k$ and we develop the iteration to end up with
\begin{equation} \label{ruggae}
\sup_{Q_{\sigma_j \theta, \sigma_j \rho}} u \leq \gamma 2^{\tilde{\gamma}j} \theta^{\frac{-2N}{N(\bar{p}-2)+2\bar{p}}} \bigg( \int \int_{Q_{\sigma_{j+1}\theta,\sigma_{j+1}\rho}} u_+^2\, dx  \bigg)^{\frac{\bar{p}}{N(\bar{p}-2)+2\bar{p}}} + \bigg(\frac{\theta^2}{\rho^p} \bigg)^{\frac{1}{2-p}}.
\end{equation} \noindent 

 \noindent 
{\small{\it STEP 2. SECOND ITERATION (condition \eqref{chi} enters)}}
\vskip0.2cm 

\noindent At this stage we would like to get an estimate with whatever power $1\leq l\leq 2$ in the integral on the right, so that we collect $\sup_{Q_{\theta,\rho}} u=M$ to get
\[
\sup_{Q_{\sigma_j \theta, \sigma_j \rho}} u \leq \gamma M_{j+1}^{\frac{(2-l)\bar{p}}{N(\bar{p}-2)+2\bar{p}}} 2^{\tilde{\gamma}j} \theta^{\frac{-2N}{N(\bar{p}-2)+2\bar{p}}} \bigg( \int \int_{Q_{\sigma_{j+1}\theta,\sigma_{j+1}\rho}} u_+^l\, dx  \bigg)^{\frac{\bar{p}}{N(\bar{p}-2)+2\bar{p}}} + \bigg(\frac{\theta^2}{\rho^p} \bigg)^{\frac{1}{2-p}}.
\] Finally we use Young's inequality with $p= \frac{N(\bar{p}-2)+2\bar{p}}{(2-l)\bar{p}}$ and $p'= (1-1/p)^{-1}= \frac{N(\bar{p}-2)+2\bar{p}}{N(\bar{p}-2)+l\bar{p}}$ to get the inequality
\begin{equation} \label{raggae-got-soul}
\sup_{Q_{\sigma_j \theta, \sigma_j \rho}} u \leq \epsilon M_{j+1} +\gamma \epsilon^{-\gamma} 2^{\tilde{\gamma} j} \theta^{\frac{-2N}{N(\bar{p}-2)+l\bar{p}}} \bigg( \int \int_{Q_{\sigma_{j+1}\theta,\sigma_{j+1}\rho}} u_+^l\, dx  \bigg)^{\frac{\bar{p}}{N(\bar{p}-2)+l\bar{p}}} + \bigg(\frac{\theta^2}{\rho^p} \bigg)^{\frac{1}{2-p}}.
\end{equation}
\noindent Now we perform a second iteration on $\sigma_j= \sum_{i=1}^j 2^{-(i+1)}$. Clearly $\sigma_0= 1/2$ and $\sigma_{\infty}=1$ and the polydiscs $ Q_{\sigma_j}$ increase up to $Q_{\theta,\rho}$. With these stipulations previous formula \eqref{raggae-got-soul} can be written as
\begin{equation*}
    \begin{aligned}
Y_{j} \leq& \epsilon Y_{j+1}+ \gamma 2^{\tilde{\gamma} j} \bigg{\{} \epsilon^{-\gamma}  \theta^{\frac{-2N}{N(\bar{p}-2)+l\bar{p}}} \bigg( \int \int_{Q_{\sigma_{j+1}\theta,\sigma_{j+1}\rho}} u_+^l\, dx  \bigg)^{\frac{\bar{p}}{N(\bar{p}-2)+l\bar{p}}} + \bigg(\frac{\theta^2}{\rho^p} \bigg)^{\frac{1}{2-p}}\bigg{\}}\\
&=\epsilon Y_{j+1}+b^j I, \quad \quad b=2^{\tilde{\gamma}}.
    \end{aligned}
\end{equation*} \noindent We iterate as in \eqref{sqn} and we obtain the inequality
\begin{equation}\label{raggae4}
\begin{aligned}
\sup_{Q_{\theta/2},\rho/2} u \leq \gamma \theta^{\frac{-2N}{N(\bar{p}-2)+l\bar{p}}} \bigg( \int \int_{Q_{\theta,\rho}} u_+^l\, dx  \bigg)^{\frac{\bar{p}}{N(\bar{p}-2)+l\bar{p}}} + \bigg(\frac{\theta^2}{\rho^p} \bigg)^{\frac{1}{2-p}}.
\end{aligned}
\end{equation} \noindent 
If we set $ \bar{p}(N-s) /p=A$, then 
\[
\frac{(N-s) \bar{p}}{p}=A= -\bigg(\frac{s\bar{p}-2N}{2}\bigg),
\] and inequality \eqref{raggae4} writes

\begin{equation}\label{raggae5}
\begin{aligned}
\sup_{Q_{\theta/2},\rho/2} u &\leq \gamma \theta^{\frac{-2N}{N(\bar{p}-2)+l\bar{p}}} (\theta^s \rho^{N-s})^{\frac{\bar{p}}{N(\bar{p}-2)+l\bar{p}}} \bigg( \dashint  \dashint_{Q_{\theta,\rho}} u_+^l\, dx  \bigg)^{\frac{\bar{p}}{N(\bar{p}-2)+l\bar{p}}} +\gamma  \bigg(\frac{\theta^2}{\rho^p} \bigg)^{\frac{1}{2-p}} \\
 \leq \gamma & \bigg( \frac{\rho^p}{\theta^2} \bigg)^{\frac{A}{N(\bar{p}-2)+\bar{p}l}}  \bigg( \dashint  \dashint_{Q_{\theta,\rho}} u_+^l\, dx  \bigg)^{\frac{\bar{p}}{N(\bar{p}-2)+l\bar{p}}} +\gamma \bigg(\frac{\theta^2}{\rho^p} \bigg)^{\frac{1}{2-p}}.
\end{aligned}
\end{equation}\noindent
\end{proof}

% \subsection{Consequences of Lemma \ref{DG-IT}}
% \begin{itemize}
%     \item When $l=1$ we have, by previous remark that $\frac{(N-s)}{\chi}=[\frac{(N-s) \bar{p}}{p[N(\bar{p}-2)+\bar{p}}]= \frac{A}{N(\bar{p}-2)+\bar{p}}$ and therefore
%   \begin{equation}\label{raggae6}
% \sup_{Q_{\theta/2},\rho/2} u 
%  \leq \gamma  \bigg( \frac{\rho^p}{\theta^2} \bigg)^{\frac{(N-s)}{\chi}}  \bigg( \dashint  \dashint_{Q_{\theta,\rho}} u_+\, dx  \bigg)^{\frac{p}{\chi}} +\gamma \bigg(\frac{\theta^2}{\rho^p} \bigg)^{\frac{1}{2-p}}
% \end{equation} \noindent as in page 20 of the old manuscript.
% \item When $l=2$ we may just stop at \eqref{ruggae} without the second iteration to get 
% \begin{equation}
%     \begin{aligned}
%     \sup_{Q_{\theta/2, \rho/2}} u &\leq \gamma \theta^{\frac{-2N}{N(\bar{p}-2)+2\bar{p}}} \bigg( \int \int_{Q_{\sigma_{j+1}\theta,\sigma_{j+1}\rho}} u_+^2\, dx  \bigg)^{\frac{\bar{p}}{N(\bar{p}-2)+2\bar{p}}} + \bigg(\frac{\theta^2}{\rho^p} \bigg)^{\frac{1}{2-p}}\\
%     &= \gamma \bigg( \frac{\rho^p}{\theta^2} \bigg)^{\frac{A}{N(\bar{p}-2)+2\bar{p}}}  \bigg( \dashint  \dashint_{Q_{\theta,\rho}} u_+^2\, dx  \bigg)^{\frac{\bar{p}}{N(\bar{p}-2)+2\bar{p}}} +\gamma \bigg(\frac{\theta^2}{\rho^p} \bigg)^{\frac{1}{2-p}}
%     \end{aligned}
% \end{equation}
% \end{itemize} \noindent 
% and therefore to get \eqref{DG} and formula (5.9) page 24 of the old manuscript.

\subsection*{Acknowledgements}
Simone Ciani and Vincenzo Vespri are members of GNAMPA (INDAM). We are grateful to Laura Baldelli, Naian Liao, Sunra N.J. Mosconi for interesting conversations on the subject. Moreover we thank Paolo Marcellini for his advice on the origins of the anisotropic problem and we acknowledge the valuable suggestions of the referees, that have improved the quality of the present paper.


\small


\begin{thebibliography}{99}
%%%%%%%%%%%%%%%%%%%%%%%%%%%%%%%%%
\bibitem{AS} S.N. Antontsev, J.I. Díaz, S. Shmarev, {\it Energy Methods for Free Boundary Problems: Applications to Nonlinear PDEs and Fluid Mechanics.} Progress in Nonlinear Differential Equations and Their Applications, Vol 48. Appl. Mech. Rev., {\bf 55}(4), 2002.

% %%%%%%%%%%%%%%%%%%
% \bibitem{Zecca} G. di Blasio, F. Feo, G. Zecca, {\it Regularity results for local solutions to some anisotropic elliptic equations.} arXiv preprint arXiv:2011.13412. (2020)


%%%%%%%%%%%%%%%%%%%%%%%%%%%%%%%%%%%%%%%%%%%%%%%%%%%%%%
\bibitem{Marcellini-Boccardo} L. Boccardo, P. Marcellini, { \it $L^{\infty}$-Regularity for Variational Problems with Sharp Non Standard Growth Conditions.} Bollettino della Unione Matematica Italiana, {\bf{7}}(4-A), (1990), 219-226.
%%%%%%%%%%%%%%%%%%%%%%%%%%%%%%%%%%%%%%%%%%%%%%%5    
\bibitem{Brasco2} P. Bousquet, L. Brasco, C. Leone, A. Verde, {\it Gradient estimates for an orthotropic nonlinear diffusion equation}. Advances in Calculus of Variations, (2021). https://doi.org/10.1515/acv-2021-0052




%%%%%%%%%%%%%%%%%%%%%%%%%%%%%%




%%%%%%%%%%%%%%%%%%%%%%%%%%%%
\bibitem{Brasco} P. Bousquet, L. Brasco, {\it Lipschitz regularity for orthotropic functionals with nonstandard growth conditions.} Rev. Mat. Iberoamericana, (2020), 3{\bf 6}(7), 1989-2032.
%%%%%%%%%%%%%%%%%%%%%%%%%%%%%%%%%%%%%%%%%%%%%%
\bibitem{Ciani} S. Ciani, S.J. Mosconi, V. Vespri, {\it Parabolic Harnack estimates for anisotropic slow diffusion.} Accepted by Journal d'Analyse Math\'ematique. https://arxiv.org/pdf/2012.09685.pdf (2020).
%%%%%%%%%%%%%%%%%%%%%%%
\bibitem{DB-Chen} Y. Chen, E. DiBenedetto, {\it H\"older estimates of solutions of singular parabolic equations with measurable coefficients.} Archive for Rational Mechanics and Analysis, {\bf{118}}(3), (1992), 257-271.



%%%%%%%%%%%%%%%%%%%%%%%%%%%%%%%%%%%%%
\bibitem{Io1} S. Ciani, V. Vespri, {\it A new short proof of regularity for local weak solutions for a certain class of singular parabolic equations.} Rendiconti di Matematica e delle sue Applicazioni, {\bf{41}}, (2020), 251-264.
%%%%%%%%%%%%%%%%%%%%%%%%%%%%%%%%%%%%%%
% \bibitem{ME} S. Ciani, V. Vespri, {\it On H\"older continuity and equivalent formulation of intrinsic Harnack estimates for an anisotropic parabolic degenerate prototype equation.} Constructive Mathematical Analysis, {\bf 4}(1), (2020), 93-103.



%%%%%%%%%%%%%%%%%%%%%%%%%%%%%%
% \bibitem{DG} E. De Giorgi, Sulla Differenziabilit\`a e l'Analiticit\`a degli Integrali Multipli Regolari,
% {\it Mem. Accad. Sci. Torino Cl. Sci. Fis. Mat. Natur.}, {\bf3}(3), (1957), 25--43.
%%%%%%%%%%%%%%%%%%%%%%%%%%%%%%

%%%%%%%%%%%%%%%%%%%%%%%%%%%%%%
% \bibitem{DB86} E. DiBenedetto, {\it On the local behavior of solutions of degenerate parabolic
% equations with measurable coefficients}, Ann. Sc. Norm. Sup. Pisa Cl. Sc. (4), {\bf13}(3), (1986),
% 487--535.
%%%%%%%%%%%%%%%%%%%%%%%%%%%%%%%%%%%%%%%%%%%%%%%%%%%%





%%%%%%%%%%%%%%%%%%%%%%%%%%%%%%%%%%%%%%%%%%%%%%

%%%%%%%%%%%%%%%%%%%%%%%%%%%%%%%%%%%%%%%
\bibitem{CU1} G. Cupini, P. Marcellini, E.  Mascolo, {\it Local boundedness of minimizers with limit growth conditions.} Journal of Optimization Theory and Applications, {\bf{166}}(1), (2015), 1-22.
%%%%%%%%%%%%%%%%%%%%%%%%%%%%%%
 \bibitem{CU2} G. Cupini, P. Marcellini, E. Mascolo, {\it Regularity of minimizers under limit growth conditions.} Nonlinear Analysis: Theory, Methods $\&$ Applications, {\bf{153}}, (2017), 294-310.
%%%%%%%%%%%%%%%%%%%%%%%%%%%%%%%%%%%%%





\bibitem{DB} E. DiBenedetto, {\it Degenerate Parabolic 
Equations}. Springer-Verlag, New York, 1993. 
%%%%%%%%%%%%%%%%%%%%%%%%%%%%%%
\bibitem{DB2} E. DiBenedetto, {\it Partial Differential Equations}. Second Edition, Birkh\"auser, Boston, 2009.
%%%%%%%%%%%%%%%%%%%%%%%%%%%%%%

%%%%%%%%%%%%%%%%%%%%%%%%%%%%%%%%%%%%%
% \bibitem{DBGV-acta} E. DiBenedetto, U. Gianazza, V. Vespri, 
% {\it Harnack estimates for quasi-linear degenerate parabolic differential equations}, 
% Acta Math., {\bf 200}(2), (2008), 181--209. 
%%%%%%%%%%%%%%%%%%%%%%%%%%%%%%
\bibitem{DBGV-pams} E. DiBenedetto, U. Gianazza and V. Vespri, 
{\it A new approach to the expansion of positivity set of non-negative 
solutions to certain singular parabolic partial differential equations}. Proc. Amer. Math. Soc., {\bf 138}(10), (2010), 3521-3529.
%%%%%%%%%%%%%%%%%%%%%%%%%%%%%%
\bibitem{DGV-Annali} E. DiBenedetto, U. Gianazza, V. Vespri, {\it Forward, backward and elliptic Harnack inequalities for non-negative solutions to certain singular parabolic partial differential equations.} Annali della Scuola Normale Superiore di Pisa-Classe di Scienze,  {\bf 9}(5), (2010), 385-422.

%%%%%%%%%%%%%%%%%%%%%%%%%%%%%%%%%%%%%%%%%%%%%%%%%%%%
% \bibitem{DBGV-mono} E. DiBenedetto, U. Gianazza, V. Vespri, 
% {\it Harnack's Inequality for Degenerate and Singular Parabolic 
% Equations}, Springer Monographs in Mathematics, Springer-Verlag, 
% New York, 2012.
%%%%%%%%%%%%%%%%%%%%%%%%%%%%%%
\bibitem{DBGV-aniso} E. DiBenedetto, U. Gianazza, V. Vespri, 
{\it Remarks on local boundedness and local H\"older continuity 
of local weak solutions to anisotropic $p$-Laplacian type equations.}  J. Elliptic Parabol. Equ., {\bf2}(1-2), (2016), 157-169.
%%%%%%%%%%%%%%%%%%%%%%%%%%%%%%
% \bibitem{DiBe} E. DiBenedetto, {\it Degenerate Parabolic Equations}, Springer-Verlag New York, 1993.
%\bibitem{Urb} J. M. Urbano, {\it The Method of Intrinsic Scaling}, Lecture Notes in Mathematics 1930, Springer-Verlag,
%2008.
%%%%%%%%%%%%%%%%%%%%%%%%%%%%%%
% \bibitem{DMV2} F. G, D\"uzg\"un, P. Marcellini, V. Vespri,
% {\it Space expansion for a solution of an anisotropic p-Laplacian equation by using a parabolic approach},
% Riv. Math. Univ. Parma, {\bf5}, (2014), 93--111.
%%%%%%%%%%%%%%%%%%%%%%%%%%%%%%%%%%
% \bibitem{Eleuteri-Marcellini-Mascolo} M. Eleuteri, P. Marcellini and E. Mascolo, {\it Regularity for scalar integrals without structure conditions.} Advances in Calculus of Variations 2018.
%%%%%%%%%%%%%%%%%%%%%%%%%%%%%%%%%%%%%%
\bibitem{Vazquez} F. Feo, J.L. V\'azquez, B. Volzone, {\it Anisotropic $p$-Laplacian Evolution of Fast Diffusion Type.} Advanced Nonlinear Studies, {\bf 21}(3), (2021), 523-555.


%%%%%%%%%%%%%%%%%%%%%%%%%%%%%%
\bibitem{FS} N. Fusco, C. Sbordone, {\it Local boundedness of minimizers in a limit case}. 
Manuscripta Math., {\bf69}(1), (1990), 19-25.
%%%%%%%%%%%%%%%%%%%%%%%%%%%%%%
\bibitem{Giaq}M. Giaquinta, 
{\it Growth conditions and regularity, a counterexample}. Manuscripta Math., {\bf59}(2), (1987), 245-248.
%%%%%%%%%%%%%%%%%%%%%%%%%%%%%%
\bibitem{Schmeiser}  J. Haskovek, C. Schmeiser, {\it A note on the anisotropic generalizations of the Sobolev and Morrey embedding theorems}. Monash Math., {\bf158},  (2009), 71-79.

%%%%%%%%%%%%%%%%%%%%%%%%%%%%%%
\bibitem{Kolodii} I.M. Kolodii, 
{\it The boundedness of generalized solutions of elliptic differential equations}. Vestnik Moskov. Univ. Ser. I Mat. Meh., {\bf25}(5), (1970), 44-52 (Russian). 
English transl.: Moscow Univ. Math. Bull. {\bf25}(5),
(1970), 31-37.%%%%%%%%%%%%%%%%%%%%%%%%%%%%%%%%%%%%%%%%%%%%%
\bibitem{Korolev} A.G. Korolev, {\it Boundedness of generalized solutions of elliptic differential equations.} Russian Math. Surveys, {\bf{38}}, (1983), 186–187.
%%%%%%%%%%%%%%%%%%%%%%%%%%%%%%%%%%%%%%%%%%%%%
\bibitem{KS} N.V. Krylov, M.V. Safonov, {\it A certain property of solutions of parabolic equations with measurable coefficients.} Izvestiya Rossiiskoi Akademii Nauk. Seriya Matematicheskaya {\bf 44}(1), (1980), 161-175.


%%%%%%%%%%%%%%%%%%%%%%%%%%%%%%%%%%%%%%%%%%%%%%
\bibitem{KK} S.N. Kruzhkov,  I.M. Kolodii, {\it On the theory of embedding of anisotropic Sobolev spaces.} Russian Mathematical Surveys, (1983), page 188.
%%%%%%%%%%%%%%%%%%%%%%%%%%%%%%
% \bibitem{GMV1}F. G. D\"uzg\"un, P. Marcellini, V. Vespri, {\it An alternative approach to H\"older continuity of solutions to some elliptic equations},
% Nonlinear Anal., {\bf94} (2014), 133--141.
%%%%%%%%%%%%%%%%%%%%%%%%%%%%%%
%\bibitem{DiBUrbVes} E. DiBenedetto, J. M. Urbano, V. Vespri, 
%{\it Current Issues on Singular and Degenerate Evolution Equations}, 
% Evolutionary Equations,
%Vol. I, Handb. Differ. Equ., North-Holland, Amsterdam, (2004), 169--286.
%%%%%%%%%%%%%%%%%%%%%%%%%%%%%%
% \bibitem{LU} O. A. Ladyzhenskaya, N. N. Ural'tseva, {\it Linear and Quasilinear Elliptic Equations.} Academic Press, New York, 1968.
%%%%%%%%%%%%%%%%%%%%%%%%%%%%%%%%%%%%%%%%%%%%%%%%%%%%%%%%
% \bibitem{Landis} E. M. Landis,{ \it Second order equations of elliptic and parabolic type.} American Mathematical Society, Transl. M. Monograph, {\bf{171}}, 1997.


%%%%%%%%%%%%%%%%%%%%
\bibitem{Liao} N. Liao, I.I. Skrypnik, V. Vespri, {\it Local regularity for an anisotropic elliptic equation.} Calculus of Variations and Partial Differential Equations, {\bf{59}}(4), (2020), 1-31.
%%%%%%%%%%%%%%%%%%%%%%%%%%%%%%
\bibitem{Lions}
J.L. Lions, {\it Quelques m\'ethodes de r\'esolution des probl\`emes aux limites non
lin\'eaires.} Dunod, Paris, 1969.
%%%%%%%%%%%%%%%%%%%%%%%%%%%%%%%%%%%%%%%%%%
\bibitem{LS} V. Liskevich, I.I. Skrypnik, {\it H\"older continuity of
 solutions to an anisotropic elliptic equation}. Nonlinear Anal., {\bf71}(5-6), (2009), 1699-1708. 
%%%%%%%%%%%%%%%%%%%%%%%%%%%%%%
\bibitem{Marce0} P. Marcellini, {\it Un exemple de solution discontinue 
d'un problème variationnel dans le cas scalaire.} Ist. Mat. U. Dini, Firenze, 1987-88.
%%%%%%%%%%%%%%%%%%%%%%%%%%%%%%
% \bibitem{Marce1} P. Marcellini,
% {\it Regularity of minimizers of integrals of the calculus of variations with nonstandard growth conditions},
% Arch. Rational Mech. Anal., {\bf105}(3), (1989), 267--284.
% %%%%%%%%%%%%%%%%%%%%%%%%%%%%%%
% \bibitem{Marce2} P. Marcellini, {\it Regularity and existence of solutions of elliptic 
% equations with $p,q$-growth conditions},
% J. Differential Equations, {\bf90}(1), (1991), 1--30. 
%%%%%%%%%%%%%%%%%%%%%%%%%%%%%%
% \bibitem{Marce3} P. Marcellini, 
% {\it Regularity for elliptic equations with general growth conditions},
% J. Differential Equations, {\bf105}(2), (1993), 296-333.%%%%%%%%%%%%%%%%%%%%%%%%%%%%%
\bibitem{Marce4} P. Marcellini, {\it Regularity under general and p, q-growth conditions.} Discrete $\&$ Continuous Dynamical Systems-Series, {\bf{13}}(7), (2020), 2009-2031.
%%%%%%%%%%%%%%%%%%%%%%%%%%%%%%%%%%%%%%%
\bibitem{Ming} G. Mingione, V. Rădulescu, {\it Recent developments in problems with nonstandard growth and nonuniform ellipticity.} Journal of Mathematical Analysis and Applications, {\bf{501}} (1), (2021), 125-197.

%%%%%%%%%%%%%%%%%%%%%%%%%%%%%%
%\bibitem{Moser61} J. Moser,
%{\it On Harnack's theorem for elliptic differential equations}, 
%Comm. Pure Appl. Math., {\bf14}, (1961), 577--591.
%%%%%%%%%%%%%%%%%%%%%%%%%%%%%%
\bibitem{Troisi} M. Troisi, {\it Teoremi di inclusione per spazi di Sobolev non isotropi}. Ricerche Mat., {\bf18},  (1969), 3-24.
%%%%%%%%%%%%%%%%%%%%%%%%%%%%%%
%%%%%%%%%%%%%%%%%%%%%%%%%%%%%%
%\bibitem{Tro2} M. Troisi, {\it Ulteriori contributi alla teoria degli spazi di Sobolev non isotropi}, 
%Ricerche Mat., {\bf20}, (1971), 90--117.
%%%%%%%%%%%%%%%%%%%%%%%%%%%%%%

\bibitem{Urbano} J.M. Urbano, {\it The method of intrinsic scaling.} Lecture Notes in Mathematics, Springer, Berlin, Heidelberg, 2008.

\bibitem{Ural}  N. N. Ural’tseva, A. B. Urdaletova, {\it The boundedness of the gradients of generalized solutions of degenerate quasilinear nonuniformly elliptic equations.} Vest. Leningr. Univ. Math {\bf 16}, (1984), 263-270.
%%%%%%%%%%%%%%%%%%%%%%%%%%%%%%%


\end{thebibliography}
\end{document}